\input amstex
%

\def\next{AMS-SEKR}\ifx\styname\next \endinput\fi
\catcode`\@=11
\def\styname{AMS-SEKR}
\def\styversion{2.0}
{\W@{}\W@{\styname.STY - Version \styversion}\W@{}}
\hyphenation{acad-e-my acad-e-mies af-ter-thought anom-aly anom-alies
an-ti-deriv-a-tive an-tin-o-my an-tin-o-mies apoth-e-o-ses apoth-e-o-sis
ap-pen-dix ar-che-typ-al as-sign-a-ble as-sist-ant-ship as-ymp-tot-ic
asyn-chro-nous at-trib-uted at-trib-ut-able bank-rupt bank-rupt-cy
bi-dif-fer-en-tial blue-print busier busiest cat-a-stroph-ic
cat-a-stroph-i-cally con-gress cross-hatched data-base de-fin-i-tive
de-riv-a-tive dis-trib-ute dri-ver dri-vers eco-nom-ics econ-o-mist
elit-ist equi-vari-ant ex-quis-ite ex-tra-or-di-nary flow-chart
for-mi-da-ble forth-right friv-o-lous ge-o-des-ic ge-o-det-ic geo-met-ric
griev-ance griev-ous griev-ous-ly hexa-dec-i-mal ho-lo-no-my ho-mo-thetic
ideals idio-syn-crasy in-fin-ite-ly in-fin-i-tes-i-mal ir-rev-o-ca-ble
key-stroke lam-en-ta-ble light-weight mal-a-prop-ism man-u-script
mar-gin-al meta-bol-ic me-tab-o-lism meta-lan-guage me-trop-o-lis
met-ro-pol-i-tan mi-nut-est mol-e-cule mono-chrome mono-pole mo-nop-oly
mono-spline mo-not-o-nous mul-ti-fac-eted mul-ti-plic-able non-euclid-ean
non-iso-mor-phic non-smooth par-a-digm par-a-bol-ic pa-rab-o-loid
pa-ram-e-trize para-mount pen-ta-gon phe-nom-e-non post-script pre-am-ble
pro-ce-dur-al pro-hib-i-tive pro-hib-i-tive-ly pseu-do-dif-fer-en-tial
pseu-do-fi-nite pseu-do-nym qua-drat-ics quad-ra-ture qua-si-smooth
qua-si-sta-tion-ary qua-si-tri-an-gu-lar quin-tes-sence quin-tes-sen-tial
re-arrange-ment rec-tan-gle ret-ri-bu-tion retro-fit retro-fit-ted
right-eous right-eous-ness ro-bot ro-bot-ics sched-ul-ing se-mes-ter
semi-def-i-nite semi-ho-mo-thet-ic set-up se-vere-ly side-step sov-er-eign
spe-cious spher-oid spher-oid-al star-tling star-tling-ly
sta-tis-tics sto-chas-tic straight-est strange-ness strat-a-gem strong-hold
sum-ma-ble symp-to-matic syn-chro-nous topo-graph-i-cal tra-vers-a-ble
tra-ver-sal tra-ver-sals treach-ery turn-around un-at-tached un-err-ing-ly
white-space wide-spread wing-spread wretch-ed wretch-ed-ly Brown-ian
Eng-lish Euler-ian Feb-ru-ary Gauss-ian Grothen-dieck Hamil-ton-ian
Her-mit-ian Jan-u-ary Japan-ese Kor-te-weg Le-gendre Lip-schitz
Lip-schitz-ian Mar-kov-ian Noe-ther-ian No-vem-ber Rie-mann-ian
Schwarz-schild Sep-tem-ber
form per-iods Uni-ver-si-ty cri-ti-sism for-ma-lism}
\Invalid@\nofrills
\Invalid@\usualspace
\newif\ifnofrills@
\def\nofrills@#1#2{\relaxnext@
  \DN@{\ifx\next\nofrills
    \nofrills@true\let#2\relax\DN@\nofrills{\nextii@}%
  \else
    \nofrills@false\def#2{#1}\let\next@\nextii@\fi
\next@}}
\def\usualspace@#1{\ifnofrills@\def\usualspace{#1}\fi}
\def\addto#1#2{\csname \expandafter\eat@\string#1@\endcsname
  \expandafter{\the\csname \expandafter\eat@\string#1@\endcsname#2}}
\newdimen\bigsize@
\def\big@#1#2{{\hbox{$\left#2\vcenter to#1\bigsize@{}%
  \right.\nulldelimiterspace\z@\m@th$}}}
\def\big{\big@\@ne}
\def\Big{\big@{1.5}}
\def\bigg{\big@\tw@}
\def\Bigg{\big@{2.5}}
\def\raggedcenter@{\leftskip\z@ plus.4\hsize \rightskip\leftskip
 \parfillskip\z@ \parindent\z@ \spaceskip.3333em \xspaceskip.5em
 \pretolerance9999\tolerance9999 \exhyphenpenalty\@M
 \hyphenpenalty\@M \let\\\linebreak}
\def\upperspecialchars{\def\ss{SS}\let\i=I\let\j=J\let\ae\AE\let\oe\OE
  \let\o\O\let\aa\AA\let\l\L}
\def\uppercasetext@#1{%
  {\spaceskip1.2\fontdimen2\the\font plus1.2\fontdimen3\the\font
   \upperspecialchars\uctext@#1$\m@th\aftergroup\eat@$}}
\def\uctext@#1$#2${\endash@#1-\endash@$#2$\uctext@}
\def\endash@#1-#2\endash@{\uppercase{#1}\if\notempty{#2}--\endash@#2\endash@\fi}
\def\runaway@#1{\DN@{#1}\ifx\envir@\next@
  \Err@{You seem to have a missing or misspelled \string\end#1 ...}%
  \let\envir@\empty\fi}
\newif\iftemp@
\def\notempty#1{TT\fi\def\test@{#1}\ifx\test@\empty\temp@false
  \else\temp@true\fi \iftemp@}
\font@\tensmc=cmcsc10
\font@\sevenex=cmex7
\font@\sevenit=cmti7
\font@\eightrm=cmr8 
\font@\sixrm=cmr6 
\font@\eighti=cmmi8     \skewchar\eighti='177 
\font@\sixi=cmmi6       \skewchar\sixi='177   
\font@\eightsy=cmsy8    \skewchar\eightsy='60 
\font@\sixsy=cmsy6      \skewchar\sixsy='60   
\font@\eightex=cmex8
\font@\eightbf=cmbx8 
\font@\sixbf=cmbx6   
\font@\eightit=cmti8 
\font@\eightsl=cmsl8 
\font@\eightsmc=cmcsc8
\font@\eighttt=cmtt8 


\loadmsam
\loadmsbm
\loadeufm
\UseAMSsymbols
\newtoks\tenpoint@
\def\tenpoint{\normalbaselineskip12\p@
 \abovedisplayskip12\p@ plus3\p@ minus9\p@
 \belowdisplayskip\abovedisplayskip
 \abovedisplayshortskip\z@ plus3\p@
 \belowdisplayshortskip7\p@ plus3\p@ minus4\p@
 \textonlyfont@\rm\tenrm \textonlyfont@\it\tenit
 \textonlyfont@\sl\tensl \textonlyfont@\bf\tenbf
 \textonlyfont@\smc\tensmc \textonlyfont@\tt\tentt
 \textonlyfont@\bsmc\tenbsmc
 \ifsyntax@ \def\big##1{{\hbox{$\left##1\right.$}}}%
  \let\Big\big \let\bigg\big \let\Bigg\big
 \else
  \textfont\z@=\tenrm  \scriptfont\z@=\sevenrm  \scriptscriptfont\z@=\fiverm
  \textfont\@ne=\teni  \scriptfont\@ne=\seveni  \scriptscriptfont\@ne=\fivei
  \textfont\tw@=\tensy \scriptfont\tw@=\sevensy \scriptscriptfont\tw@=\fivesy
  \textfont\thr@@=\tenex \scriptfont\thr@@=\sevenex
        \scriptscriptfont\thr@@=\sevenex
  \textfont\itfam=\tenit \scriptfont\itfam=\sevenit
        \scriptscriptfont\itfam=\sevenit
  \textfont\bffam=\tenbf \scriptfont\bffam=\sevenbf
        \scriptscriptfont\bffam=\fivebf
  \setbox\strutbox\hbox{\vrule height8.5\p@ depth3.5\p@ width\z@}%
  \setbox\strutbox@\hbox{\lower.5\normallineskiplimit\vbox{%
        \kern-\normallineskiplimit\copy\strutbox}}%
 \setbox\z@\vbox{\hbox{$($}\kern\z@}\bigsize@=1.2\ht\z@
 \fi
 \normalbaselines\rm\ex@.2326ex\jot3\ex@\the\tenpoint@}
\newtoks\eightpoint@
\def\eightpoint{\normalbaselineskip10\p@
 \abovedisplayskip10\p@ plus2.4\p@ minus7.2\p@
 \belowdisplayskip\abovedisplayskip
 \abovedisplayshortskip\z@ plus2.4\p@
 \belowdisplayshortskip5.6\p@ plus2.4\p@ minus3.2\p@
 \textonlyfont@\rm\eightrm \textonlyfont@\it\eightit
 \textonlyfont@\sl\eightsl \textonlyfont@\bf\eightbf
 \textonlyfont@\smc\eightsmc \textonlyfont@\tt\eighttt
 \textonlyfont@\bsmc\eightbsmc
 \ifsyntax@\def\big##1{{\hbox{$\left##1\right.$}}}%
  \let\Big\big \let\bigg\big \let\Bigg\big
 \else
  \textfont\z@=\eightrm \scriptfont\z@=\sixrm \scriptscriptfont\z@=\fiverm
  \textfont\@ne=\eighti \scriptfont\@ne=\sixi \scriptscriptfont\@ne=\fivei
  \textfont\tw@=\eightsy \scriptfont\tw@=\sixsy \scriptscriptfont\tw@=\fivesy
  \textfont\thr@@=\eightex \scriptfont\thr@@=\sevenex
   \scriptscriptfont\thr@@=\sevenex
  \textfont\itfam=\eightit \scriptfont\itfam=\sevenit
   \scriptscriptfont\itfam=\sevenit
  \textfont\bffam=\eightbf \scriptfont\bffam=\sixbf
   \scriptscriptfont\bffam=\fivebf
 \setbox\strutbox\hbox{\vrule height7\p@ depth3\p@ width\z@}%
 \setbox\strutbox@\hbox{\raise.5\normallineskiplimit\vbox{%
   \kern-\normallineskiplimit\copy\strutbox}}%
 \setbox\z@\vbox{\hbox{$($}\kern\z@}\bigsize@=1.2\ht\z@
 \fi
 \normalbaselines\eightrm\ex@.2326ex\jot3\ex@\the\eightpoint@}

\font@\twelverm=cmr10 scaled\magstep1
\font@\twelveit=cmti10 scaled\magstep1
\font@\twelvesl=cmsl10 scaled\magstep1
\font@\twelvesmc=cmcsc10 scaled\magstep1
\font@\twelvett=cmtt10 scaled\magstep1
\font@\twelvebf=cmbx10 scaled\magstep1
\font@\twelvei=cmmi10 scaled\magstep1
\font@\twelvesy=cmsy10 scaled\magstep1
\font@\twelveex=cmex10 scaled\magstep1
\font@\twelvemsa=msam10 scaled\magstep1
\font@\twelveeufm=eufm10 scaled\magstep1
\font@\twelvemsb=msbm10 scaled\magstep1
\newtoks\twelvepoint@
\def\twelvepoint{\normalbaselineskip15\p@
 \abovedisplayskip15\p@ plus3.6\p@ minus10.8\p@
 \belowdisplayskip\abovedisplayskip
 \abovedisplayshortskip\z@ plus3.6\p@
 \belowdisplayshortskip8.4\p@ plus3.6\p@ minus4.8\p@
 \textonlyfont@\rm\twelverm \textonlyfont@\it\twelveit
 \textonlyfont@\sl\twelvesl \textonlyfont@\bf\twelvebf
 \textonlyfont@\smc\twelvesmc \textonlyfont@\tt\twelvett
 \textonlyfont@\bsmc\twelvebsmc
 \ifsyntax@ \def\big##1{{\hbox{$\left##1\right.$}}}%
  \let\Big\big \let\bigg\big \let\Bigg\big
 \else
  \textfont\z@=\twelverm  \scriptfont\z@=\tenrm  \scriptscriptfont\z@=\sevenrm
  \textfont\@ne=\twelvei  \scriptfont\@ne=\teni  \scriptscriptfont\@ne=\seveni
  \textfont\tw@=\twelvesy \scriptfont\tw@=\tensy \scriptscriptfont\tw@=\sevensy
  \textfont\thr@@=\twelveex \scriptfont\thr@@=\tenex
        \scriptscriptfont\thr@@=\tenex
  \textfont\itfam=\twelveit \scriptfont\itfam=\tenit
        \scriptscriptfont\itfam=\tenit
  \textfont\bffam=\twelvebf \scriptfont\bffam=\tenbf
        \scriptscriptfont\bffam=\sevenbf
  \setbox\strutbox\hbox{\vrule height10.2\p@ depth4.2\p@ width\z@}%
  \setbox\strutbox@\hbox{\lower.6\normallineskiplimit\vbox{%
        \kern-\normallineskiplimit\copy\strutbox}}%
 \setbox\z@\vbox{\hbox{$($}\kern\z@}\bigsize@=1.4\ht\z@
 \fi
 \normalbaselines\rm\ex@.2326ex\jot3.6\ex@\the\twelvepoint@}

\def\headfonts{\twelvepoint\bf}

\font@\fourteenrm=cmr10 scaled\magstep2
\font@\fourteenit=cmti10 scaled\magstep2
\font@\fourteensl=cmsl10 scaled\magstep2
\font@\fourteensmc=cmcsc10 scaled\magstep2
\font@\fourteentt=cmtt10 scaled\magstep2
\font@\fourteenbf=cmbx10 scaled\magstep2
\font@\fourteeni=cmmi10 scaled\magstep2
\font@\fourteensy=cmsy10 scaled\magstep2
\font@\fourteenex=cmex10 scaled\magstep2
\font@\fourteenmsa=msam10 scaled\magstep2
\font@\fourteeneufm=eufm10 scaled\magstep2
\font@\fourteenmsb=msbm10 scaled\magstep2
\newtoks\fourteenpoint@
\def\fourteenpoint{\normalbaselineskip15\p@
 \abovedisplayskip18\p@ plus4.3\p@ minus12.9\p@
 \belowdisplayskip\abovedisplayskip
 \abovedisplayshortskip\z@ plus4.3\p@
 \belowdisplayshortskip10.1\p@ plus4.3\p@ minus5.8\p@
 \textonlyfont@\rm\fourteenrm \textonlyfont@\it\fourteenit
 \textonlyfont@\sl\fourteensl \textonlyfont@\bf\fourteenbf
 \textonlyfont@\smc\fourteensmc \textonlyfont@\tt\fourteentt
 \textonlyfont@\bsmc\fourteenbsmc
 \ifsyntax@ \def\big##1{{\hbox{$\left##1\right.$}}}%
  \let\Big\big \let\bigg\big \let\Bigg\big
 \else
  \textfont\z@=\fourteenrm  \scriptfont\z@=\twelverm  \scriptscriptfont\z@=\tenrm
  \textfont\@ne=\fourteeni  \scriptfont\@ne=\twelvei  \scriptscriptfont\@ne=\teni
  \textfont\tw@=\fourteensy \scriptfont\tw@=\twelvesy \scriptscriptfont\tw@=\tensy
  \textfont\thr@@=\fourteenex \scriptfont\thr@@=\twelveex
        \scriptscriptfont\thr@@=\twelveex
  \textfont\itfam=\fourteenit \scriptfont\itfam=\twelveit
        \scriptscriptfont\itfam=\twelveit
  \textfont\bffam=\fourteenbf \scriptfont\bffam=\twelvebf
        \scriptscriptfont\bffam=\tenbf
  \setbox\strutbox\hbox{\vrule height12.2\p@ depth5\p@ width\z@}%
  \setbox\strutbox@\hbox{\lower.72\normallineskiplimit\vbox{%
        \kern-\normallineskiplimit\copy\strutbox}}%
 \setbox\z@\vbox{\hbox{$($}\kern\z@}\bigsize@=1.7\ht\z@
 \fi
 \normalbaselines\rm\ex@.2326ex\jot4.3\ex@\the\fourteenpoint@}

\def\chapheadfonts{\fourteenpoint\bf}

\font@\seventeenrm=cmr10 scaled\magstep3
\font@\seventeenit=cmti10 scaled\magstep3
\font@\seventeensl=cmsl10 scaled\magstep3
\font@\seventeensmc=cmcsc10 scaled\magstep3
\font@\seventeentt=cmtt10 scaled\magstep3
\font@\seventeenbf=cmbx10 scaled\magstep3
\font@\seventeeni=cmmi10 scaled\magstep3
\font@\seventeensy=cmsy10 scaled\magstep3
\font@\seventeenex=cmex10 scaled\magstep3
\font@\seventeenmsa=msam10 scaled\magstep3
\font@\seventeeneufm=eufm10 scaled\magstep3
\font@\seventeenmsb=msbm10 scaled\magstep3
\newtoks\seventeenpoint@
\def\seventeenpoint{\normalbaselineskip18\p@
 \abovedisplayskip21.6\p@ plus5.2\p@ minus15.4\p@
 \belowdisplayskip\abovedisplayskip
 \abovedisplayshortskip\z@ plus5.2\p@
 \belowdisplayshortskip12.1\p@ plus5.2\p@ minus7\p@
 \textonlyfont@\rm\seventeenrm \textonlyfont@\it\seventeenit
 \textonlyfont@\sl\seventeensl \textonlyfont@\bf\seventeenbf
 \textonlyfont@\smc\seventeensmc \textonlyfont@\tt\seventeentt
 \textonlyfont@\bsmc\seventeenbsmc
 \ifsyntax@ \def\big##1{{\hbox{$\left##1\right.$}}}%
  \let\Big\big \let\bigg\big \let\Bigg\big
 \else
  \textfont\z@=\seventeenrm  \scriptfont\z@=\fourteenrm  \scriptscriptfont\z@=\twelverm
  \textfont\@ne=\seventeeni  \scriptfont\@ne=\fourteeni  \scriptscriptfont\@ne=\twelvei
  \textfont\tw@=\seventeensy \scriptfont\tw@=\fourteensy \scriptscriptfont\tw@=\twelvesy
  \textfont\thr@@=\seventeenex \scriptfont\thr@@=\fourteenex
        \scriptscriptfont\thr@@=\fourteenex
  \textfont\itfam=\seventeenit \scriptfont\itfam=\fourteenit
        \scriptscriptfont\itfam=\fourteenit
  \textfont\bffam=\seventeenbf \scriptfont\bffam=\fourteenbf
        \scriptscriptfont\bffam=\twelvebf
  \setbox\strutbox\hbox{\vrule height14.6\p@ depth6\p@ width\z@}%
  \setbox\strutbox@\hbox{\lower.86\normallineskiplimit\vbox{%
        \kern-\normallineskiplimit\copy\strutbox}}%
 \setbox\z@\vbox{\hbox{$($}\kern\z@}\bigsize@=2\ht\z@
 \fi
 \normalbaselines\rm\ex@.2326ex\jot5.2\ex@\the\seventeenpoint@}

\font@\rrrrrm=cmr10 scaled\magstep4
\font@\bigtitlefont=cmbx10 scaled\magstep4

\parindent1pc
\normallineskiplimit\p@
\newdimen\indenti \indenti=2pc
\def\pageheight#1{\vsize#1}
\def\pagewidth#1{\hsize#1%
   \captionwidth@\hsize \advance\captionwidth@-2\indenti}
\pagewidth{30pc} \pageheight{47pc}
\def\topmatter{%
 \ifx\undefined\msafam
 \else\font@\eightmsa=msam8 \font@\sixmsa=msam6
   \ifsyntax@\else \addto\tenpoint{\textfont\msafam=\tenmsa
              \scriptfont\msafam=\sevenmsa \scriptscriptfont\msafam=\fivemsa}%
     \addto\eightpoint{\textfont\msafam=\eightmsa \scriptfont\msafam=\sixmsa
              \scriptscriptfont\msafam=\fivemsa}%
   \fi
 \fi
 \ifx\undefined\msbfam
 \else\font@\eightmsb=msbm8 \font@\sixmsb=msbm6
   \ifsyntax@\else \addto\tenpoint{\textfont\msbfam=\tenmsb
         \scriptfont\msbfam=\sevenmsb \scriptscriptfont\msbfam=\fivemsb}%
     \addto\eightpoint{\textfont\msbfam=\eightmsb \scriptfont\msbfam=\sixmsb
         \scriptscriptfont\msbfam=\fivemsb}%
   \fi
 \fi
 \ifx\undefined\eufmfam
 \else \font@\eighteufm=eufm8 \font@\sixeufm=eufm6
   \ifsyntax@\else \addto\tenpoint{\textfont\eufmfam=\teneufm
       \scriptfont\eufmfam=\seveneufm \scriptscriptfont\eufmfam=\fiveeufm}%
     \addto\eightpoint{\textfont\eufmfam=\eighteufm
       \scriptfont\eufmfam=\sixeufm \scriptscriptfont\eufmfam=\fiveeufm}%
   \fi
 \fi
 \ifx\undefined\eufbfam
 \else \font@\eighteufb=eufb8 \font@\sixeufb=eufb6
   \ifsyntax@\else \addto\tenpoint{\textfont\eufbfam=\teneufb
      \scriptfont\eufbfam=\seveneufb \scriptscriptfont\eufbfam=\fiveeufb}%
    \addto\eightpoint{\textfont\eufbfam=\eighteufb
      \scriptfont\eufbfam=\sixeufb \scriptscriptfont\eufbfam=\fiveeufb}%
   \fi
 \fi
 \ifx\undefined\eusmfam
 \else \font@\eighteusm=eusm8 \font@\sixeusm=eusm6
   \ifsyntax@\else \addto\tenpoint{\textfont\eusmfam=\teneusm
       \scriptfont\eusmfam=\seveneusm \scriptscriptfont\eusmfam=\fiveeusm}%
     \addto\eightpoint{\textfont\eusmfam=\eighteusm
       \scriptfont\eusmfam=\sixeusm \scriptscriptfont\eusmfam=\fiveeusm}%
   \fi
 \fi
 \ifx\undefined\eusbfam
 \else \font@\eighteusb=eusb8 \font@\sixeusb=eusb6
   \ifsyntax@\else \addto\tenpoint{\textfont\eusbfam=\teneusb
       \scriptfont\eusbfam=\seveneusb \scriptscriptfont\eusbfam=\fiveeusb}%
     \addto\eightpoint{\textfont\eusbfam=\eighteusb
       \scriptfont\eusbfam=\sixeusb \scriptscriptfont\eusbfam=\fiveeusb}%
   \fi
 \fi
 \ifx\undefined\eurmfam
 \else \font@\eighteurm=eurm8 \font@\sixeurm=eurm6
   \ifsyntax@\else \addto\tenpoint{\textfont\eurmfam=\teneurm
       \scriptfont\eurmfam=\seveneurm \scriptscriptfont\eurmfam=\fiveeurm}%
     \addto\eightpoint{\textfont\eurmfam=\eighteurm
       \scriptfont\eurmfam=\sixeurm \scriptscriptfont\eurmfam=\fiveeurm}%
   \fi
 \fi
 \ifx\undefined\eurbfam
 \else \font@\eighteurb=eurb8 \font@\sixeurb=eurb6
   \ifsyntax@\else \addto\tenpoint{\textfont\eurbfam=\teneurb
       \scriptfont\eurbfam=\seveneurb \scriptscriptfont\eurbfam=\fiveeurb}%
    \addto\eightpoint{\textfont\eurbfam=\eighteurb
       \scriptfont\eurbfam=\sixeurb \scriptscriptfont\eurbfam=\fiveeurb}%
   \fi
 \fi
 \ifx\undefined\cmmibfam
 \else \font@\eightcmmib=cmmib8 \font@\sixcmmib=cmmib6
   \ifsyntax@\else \addto\tenpoint{\textfont\cmmibfam=\tencmmib
       \scriptfont\cmmibfam=\sevencmmib \scriptscriptfont\cmmibfam=\fivecmmib}%
    \addto\eightpoint{\textfont\cmmibfam=\eightcmmib
       \scriptfont\cmmibfam=\sixcmmib \scriptscriptfont\cmmibfam=\fivecmmib}%
   \fi
 \fi
 \ifx\undefined\cmbsyfam
 \else \font@\eightcmbsy=cmbsy8 \font@\sixcmbsy=cmbsy6
   \ifsyntax@\else \addto\tenpoint{\textfont\cmbsyfam=\tencmbsy
      \scriptfont\cmbsyfam=\sevencmbsy \scriptscriptfont\cmbsyfam=\fivecmbsy}%
    \addto\eightpoint{\textfont\cmbsyfam=\eightcmbsy
      \scriptfont\cmbsyfam=\sixcmbsy \scriptscriptfont\cmbsyfam=\fivecmbsy}%
   \fi
 \fi
 \let\topmatter\relax}
\def\chapterno@{\uppercase\expandafter{\romannumeral\chaptercount@}}
\newcount\chaptercount@
\def\chapter{\nofrills@{\afterassignment\chapterno@
                        CHAPTER \global\chaptercount@=}\chapter@
 \DNii@##1{\leavevmode\hskip-\leftskip
   \rlap{\vbox to\z@{\vss\centerline{\eightpoint
   \chapter@##1\unskip}\baselineskip2pc\null}}\hskip\leftskip
   \nofrills@false}%
 \FN@\next@}
\newbox\titlebox@

\def\title{\nofrills@{\relax}\title@%
 \DNii@##1\endtitle{\global\setbox\titlebox@\vtop{\tenpoint\bf
 \raggedcenter@\ignorespaces
 \baselineskip1.3\baselineskip\title@{##1}\endgraf}%
 \ifmonograph@ \edef\next{\the\leftheadtoks}\ifx\next\empty
    \leftheadtext{##1}\fi
 \fi
 \edef\next{\the\rightheadtoks}\ifx\next\empty \rightheadtext{##1}\fi
 }\FN@\next@}
\newbox\authorbox@
\def\author#1\endauthor{\global\setbox\authorbox@
 \vbox{\tenpoint\smc\raggedcenter@\ignorespaces
 #1\endgraf}\relaxnext@ \edef\next{\the\leftheadtoks}%
 \ifx\next\empty\leftheadtext{#1}\fi}
\newbox\affilbox@
\def\affil#1\endaffil{\global\setbox\affilbox@
 \vbox{\tenpoint\raggedcenter@\ignorespaces#1\endgraf}}
\newcount\addresscount@
\addresscount@\z@
\def\address#1\endaddress{\global\advance\addresscount@\@ne
  \expandafter\gdef\csname address\number\addresscount@\endcsname
  {\vskip12\p@ minus6\p@\noindent\eightpoint\smc\ignorespaces#1\par}}
\def\email{\nofrills@{\eightpoint{\it E-mail\/}:\enspace}\email@
  \DNii@##1\endemail{%
  \expandafter\gdef\csname email\number\addresscount@\endcsname
  {\def\usualspace{{\it\enspace}}\smallskip\noindent\eightpoint\email@
  \ignorespaces##1\par}}%
 \FN@\next@}
\def\thedate@{}
\def\date#1\enddate{\gdef\thedate@{\tenpoint\ignorespaces#1\unskip}}
\def\thethanks@{}
\def\thanks#1\endthanks{\gdef\thethanks@{\eightpoint\ignorespaces#1.\unskip}}
\def\thekeywords@{}
\def\keywords{\nofrills@{{\it Key words and phrases.\enspace}}\keywords@
 \DNii@##1\endkeywords{\def\thekeywords@{\def\usualspace{{\it\enspace}}%
 \eightpoint\keywords@\ignorespaces##1\unskip.}}%
 \FN@\next@}
\def\thesubjclass@{}
\def\subjclass{\nofrills@{{\rm2010 {\it Mathematics Subject
   Classification\/}.\enspace}}\subjclass@
 \DNii@##1\endsubjclass{\def\thesubjclass@{\def\usualspace
  {{\rm\enspace}}\eightpoint\subjclass@\ignorespaces##1\unskip.}}%
 \FN@\next@}
\newbox\abstractbox@
\def\abstract{\nofrills@{{\smc Abstract.\enspace}}\abstract@
 \DNii@{\setbox\abstractbox@\vbox\bgroup\noindent$$\vbox\bgroup
  \def\envir@{abstract}\advance\hsize-2\indenti
  \usualspace@{{\enspace}}\eightpoint \noindent\abstract@\ignorespaces}%
 \FN@\next@}
\def\endabstract{\par\unskip\egroup$$\egroup}
\def\widestnumber#1#2{\begingroup\let\head\null\let\subhead\empty
   \let\subsubhead\subhead
   \ifx#1\head\global\setbox\tocheadbox@\hbox{#2.\enspace}%
   \else\ifx#1\subhead\global\setbox\tocsubheadbox@\hbox{#2.\enspace}%
   \else\ifx#1\key\bgroup\let\endrefitem@\egroup
        \key#2\endrefitem@\global\refindentwd\wd\keybox@
   \else\ifx#1\no\bgroup\let\endrefitem@\egroup
        \no#2\endrefitem@\global\refindentwd\wd\nobox@
   \else\ifx#1\page\global\setbox\pagesbox@\hbox{\quad\bf#2}%
   \else\ifx#1\item\setboxz@h{#2}\global\rosteritemwd\wdz@
        \global\advance\rosteritemwd by.5\parindent
   \else\message{\string\widestnumber is not defined for this option
   (\string#1)}%
\fi\fi\fi\fi\fi\fi\endgroup}
\newif\ifmonograph@
\def\Monograph{\monograph@true \let\headmark\rightheadtext
  \let\varindent@\indent \def\headfont@{\bf}\def\proclaimheadfont@{\smc}%
  \def\demofont@{\smc}}
\let\varindent@\indent

\newbox\tocheadbox@    \newbox\tocsubheadbox@
\newbox\tocbox@
\def\toc{\toc@{Contents}}
\def\newtocdefs{%
   \def \title##1\endtitle
       {\penaltyandskip@\z@\smallskipamount
        \hangindent\wd\tocheadbox@\noindent{\bf##1}}%
   \def \chapter##1{%
        Chapter \uppercase\expandafter{\romannumeral##1.\unskip}\enspace}%
   \def \specialhead##1\endspecialhead
       {\par\hangindent\wd\tocheadbox@ \noindent##1\par}%
   \def \head##1 ##2\endhead
       {\par\hangindent\wd\tocheadbox@ \noindent
        \if\notempty{##1}\hbox to\wd\tocheadbox@{\hfil##1\enspace}\fi
        ##2\par}%
   \def \subhead##1 ##2\endsubhead
       {\par\vskip-\parskip {\normalbaselines
        \advance\leftskip\wd\tocheadbox@
        \hangindent\wd\tocsubheadbox@ \noindent
        \if\notempty{##1}\hbox to\wd\tocsubheadbox@{##1\unskip\hfil}\fi
         ##2\par}}%
   \def \subsubhead##1 ##2\endsubsubhead
       {\par\vskip-\parskip {\normalbaselines
        \advance\leftskip\wd\tocheadbox@
        \hangindent\wd\tocsubheadbox@ \noindent
        \if\notempty{##1}\hbox to\wd\tocsubheadbox@{##1\unskip\hfil}\fi
        ##2\par}}}
\def\toc@#1{\relaxnext@
   \def\page##1%
       {\unskip\penalty0\null\hfil
        \rlap{\hbox to\wd\pagesbox@{\quad\hfil##1}}\hfilneg\penalty\@M}%
 \DN@{\ifx\next\nofrills\DN@\nofrills{\nextii@}%
      \else\DN@{\nextii@{{#1}}}\fi
      \next@}%
 \DNii@##1{%
\ifmonograph@\bgroup\else\setbox\tocbox@\vbox\bgroup
   \centerline{\headfont@\ignorespaces##1\unskip}\nobreak
   \vskip\belowheadskip \fi
   \setbox\tocheadbox@\hbox{0.\enspace}%
   \setbox\tocsubheadbox@\hbox{0.0.\enspace}%
   \leftskip\indenti \rightskip\leftskip
   \setbox\pagesbox@\hbox{\bf\quad000}\advance\rightskip\wd\pagesbox@
   \newtocdefs
 }%
 \FN@\next@}
\def\endtoc{\par\egroup}
\let\pretitle\relax
\let\preauthor\relax
\let\preaffil\relax
\let\predate\relax
\let\preabstract\relax
\let\prepaper\relax
\def\dedicatory #1\enddedicatory{\def\preabstract{{\medskip
  \eightpoint\it \raggedcenter@#1\endgraf}}}
\def\thetranslator@{}
\def\translator#1\endtranslator{\def\thetranslator@{\nobreak\medskip
 \line{\eightpoint\hfil Translated by \uppercase{#1}\qquad\qquad}\nobreak}}
\outer\def\endtopmatter{\runaway@{abstract}%
 \edef\next{\the\leftheadtoks}\ifx\next\empty
  \expandafter\leftheadtext\expandafter{\the\rightheadtoks}\fi
 \ifmonograph@\else
   \ifx\thesubjclass@\empty\else \makefootnote@{}{\thesubjclass@}\fi
   \ifx\thekeywords@\empty\else \makefootnote@{}{\thekeywords@}\fi
   \ifx\thethanks@\empty\else \makefootnote@{}{\thethanks@}\fi
 \fi
  \pretitle
  \ifmonograph@ \topskip7pc \else \topskip4pc \fi
  \box\titlebox@
  \topskip10pt
  \preauthor
  \ifvoid\authorbox@\else \vskip2.5pc plus1pc \unvbox\authorbox@\fi
  \preaffil
  \ifvoid\affilbox@\else \vskip1pc plus.5pc \unvbox\affilbox@\fi
  \predate
  \ifx\thedate@\empty\else \vskip1pc plus.5pc \line{\hfil\thedate@\hfil}\fi
  \preabstract
  \ifvoid\abstractbox@\else \vskip1.5pc plus.5pc \unvbox\abstractbox@ \fi
  \ifvoid\tocbox@\else\vskip1.5pc plus.5pc \unvbox\tocbox@\fi
  \prepaper
  \vskip2pc plus1pc
}
\def\document{\let\fontlist@\relax\let\alloclist@\relax
  \tenpoint}

\newskip\aboveheadskip       \aboveheadskip1.8\bigskipamount
\newdimen\belowheadskip      \belowheadskip1.8\medskipamount

\def\headfont@{\smc}
\def\penaltyandskip@#1#2{\relax\ifdim\lastskip<#2\relax\removelastskip
      \ifnum#1=\z@\else\penalty@#1\relax\fi\vskip#2%
  \else\ifnum#1=\z@\else\penalty@#1\relax\fi\fi}
\def\nobreak{\penalty\@M
  \ifvmode\def\penalty@{\let\penalty@\penalty\count@@@}%
  \everypar{\let\penalty@\penalty\everypar{}}\fi}
\let\penalty@\penalty
\def\heading#1\endheading{\head#1\endhead}

\def\specialheadfont@{\bf}
\outer\def\specialhead{\par\penaltyandskip@{-200}\aboveheadskip
  \begingroup\interlinepenalty\@M\rightskip\z@ plus\hsize \let\\\linebreak
  \specialheadfont@\noindent\ignorespaces}
\def\endspecialhead{\par\endgroup\nobreak\vskip\belowheadskip}
\let\headmark\eat@
\newskip\subheadskip       \subheadskip\medskipamount
\def\subheadfont@{\bf}
\outer\def\subhead{\nofrills@{.\enspace}\subhead@
 \DNii@##1\endsubhead{\par\penaltyandskip@{-100}\subheadskip
  \varindent@{\usualspace@{{\subheadfont@\enspace}}%
 \subheadfont@\ignorespaces##1\unskip\subhead@}\ignorespaces}%
 \FN@\next@}
\outer\def\subsubhead{\nofrills@{.\enspace}\subsubhead@
 \DNii@##1\endsubsubhead{\par\penaltyandskip@{-50}\medskipamount
      {\usualspace@{{\it\enspace}}%
  \it\ignorespaces##1\unskip\subsubhead@}\ignorespaces}%
 \FN@\next@}
\def\proclaimheadfont@{\bf}
\outer\def\proclaim{\runaway@{proclaim}\def\envir@{proclaim}%
  \nofrills@{.\enspace}\proclaim@
 \DNii@##1{\penaltyandskip@{-100}\medskipamount\varindent@
   \usualspace@{{\proclaimheadfont@\enspace}}\proclaimheadfont@
   \ignorespaces##1\unskip\proclaim@
  \sl\ignorespaces}%
 \FN@\next@}
\outer\def\endproclaim{\let\envir@\relax\par\rm
  \penaltyandskip@{55}\medskipamount}
\def\demoheadfont@{\it}
\def\demo{\runaway@{proclaim}\nofrills@{.\enspace}\demo@
     \DNii@##1{\par\penaltyandskip@\z@\medskipamount
  {\usualspace@{{\demoheadfont@\enspace}}%
  \varindent@\demoheadfont@\ignorespaces##1\unskip\demo@}\rm
  \ignorespaces}\FN@\next@}
\def\enddemo{\par\medskip}
\def\qed{\ifhmode\unskip\nobreak\fi\quad\ifmmode\square\else$\m@th\square$\fi}
\let\remark\demo
\let\endremark\enddemo
\def\definition{\runaway@{proclaim}%
  \nofrills@{.\demoheadfont@\enspace}\definition@
        \DNii@##1{\penaltyandskip@{-100}\medskipamount
        {\usualspace@{{\demoheadfont@\enspace}}%
        \varindent@\demoheadfont@\ignorespaces##1\unskip\definition@}%
        \rm \ignorespaces}\FN@\next@}


\newdimen\rosteritemwd
\newcount\rostercount@
\newif\iffirstitem@
\let\plainitem@\item
\newtoks\everypartoks@
\def\par@{\everypartoks@\expandafter{\the\everypar}\everypar{}}
\def\roster{\edef\leftskip@{\leftskip\the\leftskip}%
 \relaxnext@
 \rostercount@\z@  
 \def\item{\FN@\rosteritem@}%
 \DN@{\ifx\next\runinitem\let\next@\nextii@\else
  \let\next@\nextiii@\fi\next@}%
 \DNii@\runinitem  
  {\unskip  
   \DN@{\ifx\next[\let\next@\nextii@\else
    \ifx\next"\let\next@\nextiii@\else\let\next@\nextiv@\fi\fi\next@}%
   \DNii@[####1]{\rostercount@####1\relax
    \enspace{\rm(\number\rostercount@)}~\ignorespaces}%
   \def\nextiii@"####1"{\enspace{\rm####1}~\ignorespaces}%
   \def\nextiv@{\enspace{\rm(1)}\rostercount@\@ne~}%
   \par@\firstitem@false  
   \FN@\next@}%
 \def\nextiii@{\par\par@  
  \penalty\@m\smallskip\vskip-\parskip
  \firstitem@true}%
 \FN@\next@}
\def\rosteritem@{\iffirstitem@\firstitem@false\else\par\vskip-\parskip\fi
 \leftskip3\parindent\noindent  
 \DNii@[##1]{\rostercount@##1\relax
  \llap{\hbox to2.5\parindent{\hss\rm(\number\rostercount@)}%
   \hskip.5\parindent}\ignorespaces}%
 \def\nextiii@"##1"{%
  \llap{\hbox to2.5\parindent{\hss\rm##1}\hskip.5\parindent}\ignorespaces}%
 \def\nextiv@{\advance\rostercount@\@ne
  \llap{\hbox to2.5\parindent{\hss\rm(\number\rostercount@)}%
   \hskip.5\parindent}}%
 \ifx\next[\let\next@\nextii@\else\ifx\next"\let\next@\nextiii@\else
  \let\next@\nextiv@\fi\fi\next@}

\newif\ifnextRunin@
\def\endroster{\relaxnext@
 \par\leftskip@  
 \penalty-50 \vskip-\parskip\smallskip  
 \DN@{\ifx\next\Runinitem\let\next@\relax
  \else\nextRunin@false\let\item\plainitem@  
   \ifx\next\par 
    \DN@\par{\everypar\expandafter{\the\everypartoks@}}%
   \else  
    \DN@{\noindent\everypar\expandafter{\the\everypartoks@}}%
  \fi\fi\next@}%
 \FN@\next@}
\newcount\rosterhangafter@
\def\Runinitem#1\roster\runinitem{\relaxnext@
 \rostercount@\z@ 
 \def\item{\FN@\rosteritem@}%
 \def\runinitem@{#1}%
 \DN@{\ifx\next[\let\next\nextii@\else\ifx\next"\let\next\nextiii@
  \else\let\next\nextiv@\fi\fi\next}%
 \DNii@[##1]{\rostercount@##1\relax
  \def\item@{{\rm(\number\rostercount@)}}\nextv@}%
 \def\nextiii@"##1"{\def\item@{{\rm##1}}\nextv@}%
 \def\nextiv@{\advance\rostercount@\@ne
  \def\item@{{\rm(\number\rostercount@)}}\nextv@}%
 \def\nextv@{\setbox\z@\vbox  
  {\ifnextRunin@\noindent\fi  
  \runinitem@\unskip\enspace\item@~\par  
  \global\rosterhangafter@\prevgraf}%
  \firstitem@false  
  \ifnextRunin@\else\par\fi  
  \hangafter\rosterhangafter@\hangindent3\parindent
  \ifnextRunin@\noindent\fi  
  \runinitem@\unskip\enspace 
  \item@~\ifnextRunin@\else\par@\fi  
  \nextRunin@true\ignorespaces}%
 \FN@\next@}
\def\footmarkform@#1{$\m@th^{#1}$}
\let\thefootnotemark\footmarkform@
\def\makefootnote@#1#2{\insert\footins
 {\interlinepenalty\interfootnotelinepenalty
 \eightpoint\splittopskip\ht\strutbox\splitmaxdepth\dp\strutbox
 \floatingpenalty\@MM\leftskip\z@\rightskip\z@\spaceskip\z@\xspaceskip\z@
 \leavevmode{#1}\footstrut\ignorespaces#2\unskip\lower\dp\strutbox
 \vbox to\dp\strutbox{}}}
\newcount\footmarkcount@
\footmarkcount@\z@
\def\footnotemark{\let\@sf\empty\relaxnext@
 \ifhmode\edef\@sf{\spacefactor\the\spacefactor}\/\fi
 \DN@{\ifx[\next\let\next@\nextii@\else
  \ifx"\next\let\next@\nextiii@\else
  \let\next@\nextiv@\fi\fi\next@}%
 \DNii@[##1]{\footmarkform@{##1}\@sf}%
 \def\nextiii@"##1"{{##1}\@sf}%
 \def\nextiv@{\iffirstchoice@\global\advance\footmarkcount@\@ne\fi
  \footmarkform@{\number\footmarkcount@}\@sf}%
 \FN@\next@}
\def\footnotetext{\relaxnext@
 \DN@{\ifx[\next\let\next@\nextii@\else
  \ifx"\next\let\next@\nextiii@\else
  \let\next@\nextiv@\fi\fi\next@}%
 \DNii@[##1]##2{\makefootnote@{\footmarkform@{##1}}{##2}}%
 \def\nextiii@"##1"##2{\makefootnote@{##1}{##2}}%
 \def\nextiv@##1{\makefootnote@{\footmarkform@{\number\footmarkcount@}}{##1}}%
 \FN@\next@}
\def\footnote{\let\@sf\empty\relaxnext@
 \ifhmode\edef\@sf{\spacefactor\the\spacefactor}\/\fi
 \DN@{\ifx[\next\let\next@\nextii@\else
  \ifx"\next\let\next@\nextiii@\else
  \let\next@\nextiv@\fi\fi\next@}%
 \DNii@[##1]##2{\footnotemark[##1]\footnotetext[##1]{##2}}%
 \def\nextiii@"##1"##2{\footnotemark"##1"\footnotetext"##1"{##2}}%
 \def\nextiv@##1{\footnotemark\footnotetext{##1}}%
 \FN@\next@}
\def\adjustfootnotemark#1{\advance\footmarkcount@#1\relax}
\def\footnoterule{\kern-3\p@
  \hrule width 5pc\kern 2.6\p@} 
\def\captionfont@{\smc}
\def\topcaption#1#2\endcaption{%
  {\dimen@\hsize \advance\dimen@-\captionwidth@
   \rm\raggedcenter@ \advance\leftskip.5\dimen@ \rightskip\leftskip
  {\captionfont@#1}%
  \if\notempty{#2}.\enspace\ignorespaces#2\fi
  \endgraf}\nobreak\bigskip}
\def\botcaption#1#2\endcaption{%
  \nobreak\bigskip
  \setboxz@h{\captionfont@#1\if\notempty{#2}.\enspace\rm#2\fi}%
  {\dimen@\hsize \advance\dimen@-\captionwidth@
   \leftskip.5\dimen@ \rightskip\leftskip
   \noindent \ifdim\wdz@>\captionwidth@ 
   \else\hfil\fi 
  {\captionfont@#1}\if\notempty{#2}.\enspace\rm#2\fi\endgraf}}
\def\@ins{\par\begingroup\def\vspace##1{\vskip##1\relax}%
  \def\captionwidth##1{\captionwidth@##1\relax}%
  \setbox\z@\vbox\bgroup} 
\def\block{\RIfMIfI@\nondmatherr@\block\fi
       \else\ifvmode\vskip\abovedisplayskip\noindent\fi
        $$\def\endblock{\par\egroup$$}\fi
  \vbox\bgroup\advance\hsize-2\indenti\noindent}
\def\endblock{\par\egroup}
\def\cite#1{{\rm[{\citefont@\m@th#1}]}}
\def\citefont@{\rm}
\def\refsfont@{\eightpoint}
\outer\def\Refs{\runaway@{proclaim}%
 \relaxnext@ \DN@{\ifx\next\nofrills\DN@\nofrills{\nextii@}\else
  \DN@{\nextii@{References}}\fi\next@}%
 \DNii@##1{\penaltyandskip@{-200}\aboveheadskip
  \line{\hfil\headfont@\ignorespaces##1\unskip\hfil}\nobreak
  \vskip\belowheadskip
  \begingroup\refsfont@\sfcode`.=\@m}%
 \FN@\next@}
\def\endRefs{\par\endgroup}
\newbox\nobox@            \newbox\keybox@           \newbox\bybox@
\newbox\paperbox@         \newbox\paperinfobox@     \newbox\jourbox@
\newbox\volbox@           \newbox\issuebox@         \newbox\yrbox@
\newbox\pagesbox@         \newbox\bookbox@          \newbox\bookinfobox@
\newbox\publbox@          \newbox\publaddrbox@      \newbox\finalinfobox@
\newbox\edsbox@           \newbox\langbox@
\newif\iffirstref@        \newif\iflastref@
\newif\ifprevjour@        \newif\ifbook@            \newif\ifprevinbook@
\newif\ifquotes@          \newif\ifbookquotes@      \newif\ifpaperquotes@
\newdimen\bysamerulewd@
\setboxz@h{\refsfont@\kern3em}
\bysamerulewd@\wdz@
\newdimen\refindentwd
\setboxz@h{\refsfont@ 00. }
\refindentwd\wdz@
\outer\def\ref{\begingroup \noindent\hangindent\refindentwd
 \firstref@true \def\nofrills{\def\refkern@{\kern3sp}}%
 \ref@}
\def\ref@{\book@false \bgroup\let\endrefitem@\egroup \ignorespaces}
\def\moreref{\endrefitem@\endref@\firstref@false\ref@}%
\def\transl{\endrefitem@\endref@\firstref@false
  \book@false
  \prepunct@
  \setboxz@h\bgroup \aftergroup\unhbox\aftergroup\z@
    \def\endrefitem@{\unskip\refkern@\egroup}\ignorespaces}%
\def\emptyifempty@{\dimen@\wd\currbox@
  \advance\dimen@-\wd\z@ \advance\dimen@-.1\p@
  \ifdim\dimen@<\z@ \setbox\currbox@\copy\voidb@x \fi}
\let\refkern@\relax
\def\endrefitem@{\unskip\refkern@\egroup
  \setboxz@h{\refkern@}\emptyifempty@}\ignorespaces
\def\refdef@#1#2#3{\edef\next@{\noexpand\endrefitem@
  \let\noexpand\currbox@\csname\expandafter\eat@\string#1box@\endcsname
    \noexpand\setbox\noexpand\currbox@\hbox\bgroup}%
  \toks@\expandafter{\next@}%
  \if\notempty{#2#3}\toks@\expandafter{\the\toks@
  \def\endrefitem@{\unskip#3\refkern@\egroup
  \setboxz@h{#2#3\refkern@}\emptyifempty@}#2}\fi
  \toks@\expandafter{\the\toks@\ignorespaces}%
  \edef#1{\the\toks@}}
\refdef@\no{}{. }
\refdef@\key{[\m@th}{] }
\refdef@\by{}{}
\def\bysame{\by\hbox to\bysamerulewd@{\hrulefill}\thinspace
   \kern0sp}
\def\manyby{\message{\string\manyby is no longer necessary; \string\by
  can be used instead, starting with version 2.0 of \styname.STY}\by}
\refdef@\paper{\ifpaperquotes@``\fi\it}{}
\refdef@\paperinfo{}{}
\def\jour{\endrefitem@\let\currbox@\jourbox@
  \setbox\currbox@\hbox\bgroup
  \def\endrefitem@{\unskip\refkern@\egroup
    \setboxz@h{\refkern@}\emptyifempty@
    \ifvoid\jourbox@\else\prevjour@true\fi}%
\ignorespaces}
\refdef@\vol{\ifbook@\else\bf\fi}{}
\refdef@\issue{no. }{}
\refdef@\yr{}{}
\refdef@\pages{}{}
\def\page{\endrefitem@\def\pp@{\def\pp@{pp.~}p.~}\let\currbox@\pagesbox@
  \setbox\currbox@\hbox\bgroup\ignorespaces}
\def\pp@{pp.~}
\def\book{\endrefitem@ \let\currbox@\bookbox@
 \setbox\currbox@\hbox\bgroup\def\endrefitem@{\unskip\refkern@\egroup
  \setboxz@h{\ifbookquotes@``\fi}\emptyifempty@
  \ifvoid\bookbox@\else\book@true\fi}%
  \ifbookquotes@``\fi\it\ignorespaces}
\def\inbook{\endrefitem@
  \let\currbox@\bookbox@\setbox\currbox@\hbox\bgroup
  \def\endrefitem@{\unskip\refkern@\egroup
  \setboxz@h{\ifbookquotes@``\fi}\emptyifempty@
  \ifvoid\bookbox@\else\book@true\previnbook@true\fi}%
  \ifbookquotes@``\fi\ignorespaces}
\refdef@\eds{(}{, eds.)}
\def\ed{\endrefitem@\let\currbox@\edsbox@
 \setbox\currbox@\hbox\bgroup
 \def\endrefitem@{\unskip, ed.)\refkern@\egroup
  \setboxz@h{(, ed.)}\emptyifempty@}(\ignorespaces}
\refdef@\bookinfo{}{}
\refdef@\publ{}{}
\refdef@\publaddr{}{}
\refdef@\finalinfo{}{}
\refdef@\lang{(}{)}

\let\refdef@\relax 
\def\ppunbox@#1{\ifvoid#1\else\prepunct@\unhbox#1\fi}
\def\nocomma@#1{\ifvoid#1\else\changepunct@3\prepunct@\unhbox#1\fi}
\def\changepunct@#1{\ifnum\lastkern<3 \unkern\kern#1sp\fi}
\def\prepunct@{\count@\lastkern\unkern
  \ifnum\lastpenalty=0
    \let\penalty@\relax
  \else
    \edef\penalty@{\penalty\the\lastpenalty\relax}%
  \fi
  \unpenalty
  \let\refspace@\ \ifcase\count@,
\or;\or.\or 
  \or\let\refspace@\relax
  \else,\fi
  \ifquotes@''\quotes@false\fi \penalty@ \refspace@
}
\def\transferpenalty@#1{\dimen@\lastkern\unkern
  \ifnum\lastpenalty=0\unpenalty\let\penalty@\relax
  \else\edef\penalty@{\penalty\the\lastpenalty\relax}\unpenalty\fi
  #1\penalty@\kern\dimen@}
\def\endref{\endrefitem@\lastref@true\endref@
  \par\endgroup \prevjour@false \previnbook@false }
\def\endref@{%
\iffirstref@
  \ifvoid\nobox@\ifvoid\keybox@\indent\fi
  \else\hbox to\refindentwd{\hss\unhbox\nobox@}\fi
  \ifvoid\keybox@
  \else\ifdim\wd\keybox@>\refindentwd
         \box\keybox@
       \else\hbox to\refindentwd{\unhbox\keybox@\hfil}\fi\fi
  \kern4sp\ppunbox@\bybox@
\fi 
  \ifvoid\paperbox@
  \else\prepunct@\unhbox\paperbox@
    \ifpaperquotes@\quotes@true\fi\fi
  \ppunbox@\paperinfobox@
  \ifvoid\jourbox@
    \ifprevjour@ \nocomma@\volbox@
      \nocomma@\issuebox@
      \ifvoid\yrbox@\else\changepunct@3\prepunct@(\unhbox\yrbox@
        \transferpenalty@)\fi
      \ppunbox@\pagesbox@
    \fi 
  \else \prepunct@\unhbox\jourbox@
    \nocomma@\volbox@
    \nocomma@\issuebox@
    \ifvoid\yrbox@\else\changepunct@3\prepunct@(\unhbox\yrbox@
      \transferpenalty@)\fi
    \ppunbox@\pagesbox@
  \fi 
  \ifbook@\prepunct@\unhbox\bookbox@ \ifbookquotes@\quotes@true\fi \fi
  \nocomma@\edsbox@
  \ppunbox@\bookinfobox@
  \ifbook@\ifvoid\volbox@\else\prepunct@ vol.~\unhbox\volbox@
  \fi\fi
  \ppunbox@\publbox@ \ppunbox@\publaddrbox@
  \ifbook@ \ppunbox@\yrbox@
    \ifvoid\pagesbox@
    \else\prepunct@\pp@\unhbox\pagesbox@\fi
  \else
    \ifprevinbook@ \ppunbox@\yrbox@
      \ifvoid\pagesbox@\else\prepunct@\pp@\unhbox\pagesbox@\fi
    \fi \fi
  \ppunbox@\finalinfobox@
  \iflastref@
    \ifvoid\langbox@.\ifquotes@''\fi
    \else\changepunct@2\prepunct@\unhbox\langbox@\fi
  \else
    \ifvoid\langbox@\changepunct@1%
    \else\changepunct@3\prepunct@\unhbox\langbox@
      \changepunct@1\fi
  \fi
}
\outer\def\enddocument{%
 \runaway@{proclaim}%
\ifmonograph@ 
\else
 \nobreak
 \thetranslator@
 \count@\z@ \loop\ifnum\count@<\addresscount@\advance\count@\@ne
 \csname address\number\count@\endcsname
 \csname email\number\count@\endcsname
 \repeat
\fi
 \vfill\supereject\end}

\def\headfont@{\headfonts}
\def\proclaimheadfont@{\bf}
\def\specialheadfont@{\bf}
\def\subheadfont@{\bf}
\def\demoheadfont@{\smc}

\newif\ifThisToToc \ThisToTocfalse
\newif\iftocloaded \tocloadedfalse

\def\C@L{\noexpand\Cal}\def\B@B{\noexpand\Bbb}\def\fR@K{\noexpand\frak}
\def\S@{\noexpand\S}\def\P@P{\noexpand\"}
\def\xpar{\\}

\def\writetoc#1{\iftocloaded\ifThisToToc\begingroup\def\totoc{}
  \def\Cal{\noexpand\C@L}\def\Bbb{\noexpand\B@B}
  \def\frak{\noexpand\fR@K}\def\goth{\frak}\def\S{\noexpand\S@}
  \def\"{\noexpand\P@P}
  \def\xpar{\par\penalty100000 }\def\idx##1{##1}\def\\{\xpar}
  \edef\next@{\write\toc{\noindent#1\leaderfill\noexpand\folio\par}}%
  \next@\endgroup\global\ThisToTocfalse\fi\fi}
\def\leaderfill{\leaders\hbox to 1em{\hss.\hss}\hfill}

\newif\ifindexloaded \indexloadedfalse
\def\idx#1{\ifindexloaded\begingroup\def\ign{}\def\it{}\def\/{}%
 \def\smc{}\def\bf{}\def\tt{}%
 \def\Cal{\noexpand\C@L}\def\Bbb{\noexpand\B@B}%
 \def\frak{\noexpand\fR@K}\def\goth{\frak}\def\S{\noexpand\S@}%
  \def\"{\noexpand\P@P}%
 {\edef\next@{\write\index{#1, \noexpand\folio}}\next@}%
 \endgroup\fi{#1}}
\def\ign#1{}

\def\totoc{\global\ThisToToctrue}

\outer\def\head#1\endhead{\par\penaltyandskip@{-200}\aboveheadskip
 {\headfont@\raggedcenter@\interlinepenalty\@M
 \ignorespaces#1\endgraf}\nobreak
 \vskip\belowheadskip
 \headmark{#1}\writetoc{#1}}

\outer\def\chaphead#1\endchaphead{\par\penaltyandskip@{-200}\aboveheadskip
 {\chapheadfonts\raggedcenter@\interlinepenalty\@M
 \ignorespaces#1\endgraf}\nobreak
 \vskip3\belowheadskip
 \headmark{#1}\writetoc{#1}}

\def\folio{{\foliofont@\ifnum\pageno<\z@ \romannumeral-\pageno
 \else\number\pageno \fi}}
\newtoks\leftheadtoks
\newtoks\rightheadtoks

\def\leftheadtext{\nofrills@{\relax}\lht@
  \DNii@##1{\leftheadtoks\expandafter{\lht@{##1}}%
    \mark{\the\leftheadtoks\noexpand\else\the\rightheadtoks}
    \ifsyntax@\setboxz@h{\def\\{\unskip\space\ignorespaces}%
        \headlinefont@##1}\fi}%
  \FN@\next@}
\def\rightheadtext{\nofrills@{\relax}\rht@
  \DNii@##1{\rightheadtoks\expandafter{\rht@{##1}}%
    \mark{\the\leftheadtoks\noexpand\else\the\rightheadtoks}%
    \ifsyntax@\setboxz@h{\def\\{\unskip\space\ignorespaces}%
        \headlinefont@##1}\fi}%
  \FN@\next@}
\def\NoRunningHeads{\global\runheads@false\global\let\headmark\eat@}

\newif\iffirstpage@     \firstpage@true
\newif\ifrunheads@      \runheads@true

\newdimen\fullhsize \fullhsize=\hsize
\newdimen\fullvsize \fullvsize=\vsize
\def\fullline{\hbox to\fullhsize}

\def\pagenumbers{\gdef\folio{\folio@}}

\let\norunningheads\NoRunningHeads
\def\userunningheads{\global\runheads@true}
\norunningheads

\headline={\def\chapter#1{\chapterno@. }%
  \def\\{\unskip\space\ignorespaces}\ifrunheads@\headlinefont@
    \ifodd\pageno\rightheadline \else\leftheadline\fi
   \else\hfil\fi\ifNoRunHeadline\global\NoRunHeadlinefalse\fi}
\let\folio@\folio
\def\foliofont@{\foliofont}
\def\foliofont{\eightrm}
\def\headlinefont@{\headlinefont}
\def\headlinefont{\eightpoint\smc}
\def\leftheadline{\rlap{\folio}\hfill
   \ifNoRunHeadline\else\iftrue\topmark\fi\fi \hfill}
\def\rightheadline{\hfill\ifNoRunHeadline
   \else \expandafter\fi
  \hfill \llap{\folio}}
\footline={{\eightpoint\bottremark}%
   \ifrunheads@\else\hfil{\let\foliofont\tenrm\folio}\fi\hfil}
\def\bottremark{}
 
\newif\ifNoRunHeadline      
\def\norunninghead{\global\NoRunHeadlinetrue}
\norunninghead

\output={\output@}
%
\newif\ifoffset\offsetfalse
\output={\output@}
\def\output@{%
 \ifoffset 
  \ifodd\count0\advance\hoffset by0.5truecm
   \else\advance\hoffset by-0.5truecm\fi\fi
 \shipout\vbox{%
  \makeheadline \pagebody \makefootline }%
 \advancepageno \ifnum\outputpenalty>-\@MM\else\dosupereject\fi}

\def\indexoutput#1{%
  \ifoffset 
   \ifodd\count0\advance\hoffset by0.5truecm
    \else\advance\hoffset by-0.5truecm\fi\fi
  \shipout\vbox{\makeheadline
  \vbox to\fullvsize{\boxmaxdepth\maxdepth%
  \ifvoid\topins\else\unvbox\topins\fi%
  #1 %
  \ifvoid\footins\else 
    \vskip\skip\footins
    \footnoterule
    \unvbox\footins\fi
  \ifr@ggedbottom \kern-\dimen@ \vfil \fi}%
  \baselineskip2pc
  \makefootline}%
 \global\advance\pageno\@ne
 \ifnum\outputpenalty>-\@MM\else\dosupereject\fi}
 
 \newbox\partialpage \newdimen\halfsize \halfsize=0.5\fullhsize
 \advance\halfsize by-0.5em

 \def\begindoublecolumns{\output={\indexoutput{\unvbox255}}%
   \begingroup \def\line{\fullline}
   \output={\global\setbox\partialpage=\vbox{\unvbox255\bigskip}}\eject
   \output={\doublecolumnout}\hsize=\halfsize \vsize=2\fullvsize}
 \def\enddoublecolumns{\output={\balancecolumns}\eject
  \endgroup \pagegoal=\fullvsize%
  \output={\output@}}
\def\doublecolumnout{\splittopskip=\topskip \splitmaxdepth=\maxdepth
  \dimen@=\fullvsize \advance\dimen@ by-\ht\partialpage
  \setbox0=\vsplit255 to \dimen@ \setbox2=\vsplit255 to \dimen@
  \indexoutput{\pagesofar} \unvbox255 \penalty\outputpenalty}
\def\pagesofar{\unvbox\partialpage
  \wd0=\hsize \wd2=\hsize \hbox to\fullhsize{\box0\hfil\box2}}
\def\balancecolumns{\setbox0=\vbox{\unvbox255} \dimen@=\ht0
  \advance\dimen@ by\topskip \advance\dimen@ by-\baselineskip
  \divide\dimen@ by2 \splittopskip=\topskip
  {\vbadness=10000 \loop \global\setbox3=\copy0
    \global\setbox1=\vsplit3 to\dimen@
    \ifdim\ht3>\dimen@ \global\advance\dimen@ by1pt \repeat}
  \setbox0=\vbox to\dimen@{\unvbox1} \setbox2=\vbox to\dimen@{\unvbox3}
  \pagesofar}

\tenpoint
\catcode`\@=\active

\def\smallheadings{\let\chapheadfonts\tenpoint\let\headfonts\tenpoint}

\tenpoint
\catcode`\@=\active

\def\LL{\leavevmode\setbox0=\hbox{L}\hbox to\wd0{\hss\char'40L}}
\def\al{\alpha}

\def\la{\lambda}
\def\rh{\rho}


\def\today{\ifcase\month\or
 January\or February\or March\or April\or May\or June\or
 July\or August\or September\or October\or November\or December\fi
 \space\number\day, \number\year}

\def\({\left(}
\def\){\right)}
\def\[{\left[}
\def\]{\right]}

\def\3{\ss}
\catcode`\@=11
\def\dddot#1{\vbox{\ialign{##\crcr
      .\hskip-.5pt.\hskip-.5pt.\crcr\noalign{\kern1.5\p@\nointerlineskip}
      $\hfil\displaystyle{#1}\hfil$\crcr}}}

\newif\iftab@\tab@false
\newif\ifvtab@\vtab@false
\def\tab{\bgroup\tab@true\vtab@false\vst@bfalse\Strich@false%
   \def\\{\global\hline@@false%
     \ifhline@\global\hline@false\global\hline@@true\fi\cr}
   \edef\l@{\the\leftskip}\ialign\bgroup\hskip\l@##\hfil&&##\hfil\cr}
\def\endtab{\cr\egroup\egroup}
\def\vtab{\vtop\bgroup\vst@bfalse\vtab@true\tab@true\Strich@false%
   \bgroup\def\\{\cr}\ialign\bgroup&##\hfil\cr}
\def\endvtab{\cr\egroup\egroup\egroup}
\def\stab{\D@cke0.5pt\null 
 \bgroup\tab@true\vtab@false\vst@bfalse\Strich@true\Let@@\vspace@
 \normalbaselines\offinterlineskip
  \openup\spreadmlines@
 \edef\l@{\the\leftskip}\ialign
 \bgroup\hskip\l@##\hfil&&##\hfil\crcr}
\def\endstab{\crcr\egroup
 \egroup}
\newif\ifvst@b\vst@bfalse
\def\vstab{\D@cke0.5pt\null
 \vtop\bgroup\tab@true\vtab@false\vst@btrue\Strich@true\bgroup\Let@@\vspace@
 \normalbaselines\offinterlineskip
  \openup\spreadmlines@\bgroup}
\def\endvstab{\crcr\egroup\egroup
 \egroup\tab@false\Strich@false}

\newdimen\htstrut@
\htstrut@8.5\p@
\newdimen\htStrut@
\htStrut@12\p@
\newdimen\dpstrut@
\dpstrut@3.5\p@
\newdimen\dpStrut@
\dpStrut@3.5\p@
\def\openup{\afterassignment\@penup\dimen@=}
\def\@penup{\advance\lineskip\dimen@
  \advance\baselineskip\dimen@
  \advance\lineskiplimit\dimen@
  \divide\dimen@ by2
  \advance\htstrut@\dimen@
  \advance\htStrut@\dimen@
  \advance\dpstrut@\dimen@
  \advance\dpStrut@\dimen@}
\def\Let@@{\relax%
    \def\\{\global\hline@@false%
     \ifhline@\global\hline@false\global\hline@@true\fi\cr}%
    \iffalse}\fi}
\def\matrix{\null\,\vcenter\bgroup
 \tab@false\vtab@false\vst@bfalse\Strich@false\Let@@\vspace@
 \normalbaselines\openup\spreadmlines@\ialign
 \bgroup\hfil$\m@th##$\hfil&&\quad\hfil$\m@th##$\hfil\crcr
 \Mathstrut@\crcr\noalign{\kern-\baselineskip}}
\def\endmatrix{\crcr\Mathstrut@\crcr\noalign{\kern-\baselineskip}\egroup
 \egroup\,}
\def\smatrix{\D@cke0.5pt\null\,
 \vcenter\bgroup\tab@false\vtab@false\vst@bfalse\Strich@true\Let@@\vspace@
 \normalbaselines\offinterlineskip
  \openup\spreadmlines@\ialign
 \bgroup\hfil$\m@th##$\hfil&&\quad\hfil$\m@th##$\hfil\crcr}
\def\endsmatrix{\crcr\egroup
 \egroup\,\Strich@false}
\newdimen\D@cke
\def\Dicke#1{\global\D@cke#1}
\newtoks\tabs@\tabs@{&}
\newif\ifStrich@\Strich@false
\newif\iff@rst

\def\Stricherr@{\iftab@\ifvtab@\errmessage{\noexpand\s not allowed
     here. Use \noexpand\vstab!}%
  \else\errmessage{\noexpand\s not allowed here. Use \noexpand\stab!}%
  \fi\else\errmessage{\noexpand\s not allowed
     here. Use \noexpand\smatrix!}\fi}
\def\format{\ifvst@b\else\crcr\fi\egroup\iffalse{\fi\ifnum`}=0 \fi\format@}
\def\format@#1\\{\def\preamble@{#1}%
 \def\Str@chfehlt##1{\ifx##1\s\Stricherr@\fi\ifx##1\\\let\Next\relax%
   \else\let\Next\Str@chfehlt\fi\Next}%
 \def\c{\hfil\noexpand\ifhline@@\hbox{\vrule height\htStrut@%
   depth\dpstrut@ width\z@}\noexpand\fi%
   \ifStrich@\hbox{\vrule height\htstrut@ depth\dpstrut@ width\z@}%
   \fi\iftab@\else$\m@th\fi\the\hashtoks@\iftab@\else$\fi\hfil}%
 \def\r{\hfil\noexpand\ifhline@@\hbox{\vrule height\htStrut@%
   depth\dpstrut@ width\z@}\noexpand\fi%
   \ifStrich@\hbox{\vrule height\htstrut@ depth\dpstrut@ width\z@}%
   \fi\iftab@\else$\m@th\fi\the\hashtoks@\iftab@\else$\fi}%
 \def\l{\noexpand\ifhline@@\hbox{\vrule height\htStrut@%
   depth\dpstrut@ width\z@}\noexpand\fi%
   \ifStrich@\hbox{\vrule height\htstrut@ depth\dpstrut@ width\z@}%
   \fi\iftab@\else$\m@th\fi\the\hashtoks@\iftab@\else$\fi\hfil}%
 \def\s{\ifStrich@\ \the\tabs@\vrule width\D@cke\the\hashtoks@%
          \fi\the\tabs@\ }%
 \def\sa{\ifStrich@\vrule width\D@cke\the\hashtoks@%
            \the\tabs@\ %
            \fi}%
 \def\se{\ifStrich@\ \the\tabs@\vrule width\D@cke\the\hashtoks@\fi}%
 \def\cd{\hfil\noexpand\ifhline@@\hbox{\vrule height\htStrut@%
   depth\dpstrut@ width\z@}\noexpand\fi%
   \ifStrich@\hbox{\vrule height\htstrut@ depth\dpstrut@ width\z@}%
   \fi$\dsize\m@th\the\hashtoks@$\hfil}%
 \def\rd{\hfil\noexpand\ifhline@@\hbox{\vrule height\htStrut@%
   depth\dpstrut@ width\z@}\noexpand\fi%
   \ifStrich@\hbox{\vrule height\htstrut@ depth\dpstrut@ width\z@}%
   \fi$\dsize\m@th\the\hashtoks@$}%
 \def\ld{\noexpand\ifhline@@\hbox{\vrule height\htStrut@%
   depth\dpstrut@ width\z@}\noexpand\fi%
   \ifStrich@\hbox{\vrule height\htstrut@ depth\dpstrut@ width\z@}%
   \fi$\dsize\m@th\the\hashtoks@$\hfil}%
 \ifStrich@\else\Str@chfehlt#1\\\fi%
 \setbox\z@\hbox{\xdef\Preamble@{\preamble@}}\ifnum`{=0 \fi\iffalse}\fi
 \ialign\bgroup\span\Preamble@\crcr}
\newif\ifhline@\hline@false
\newif\ifhline@@\hline@@false
\def\hlinefor#1{\multispan@{\strip@#1 }\leaders\hrule height\D@cke\hfill%
    \global\hline@true\ignorespaces}
\def\Item "#1"{\par\noindent\hangindent2\parindent%
  \hangafter1\setbox0\hbox{\rm#1\enspace}\ifdim\wd0>2\parindent%
  \box0\else\hbox to 2\parindent{\rm#1\hfil}\fi\ignorespaces}
\def\ITEM #1"#2"{\par\noindent\hangafter1\hangindent#1%
  \setbox0\hbox{\rm#2\enspace}\ifdim\wd0>#1%
  \box0\else\hbox to 0pt{\rm#2\hss}\hskip#1\fi\ignorespaces}
\def\item"#1"{\par\noindent\hang%
  \setbox0=\hbox{\rm#1\enspace}\ifdim\wd0>\the\parindent%
  \box0\else\hbox to \parindent{\rm#1\hfil}\enspace\fi\ignorespaces}
\let\plainitem@\item
\catcode`\@=13

\magnification1200

\TagsOnRight

\catcode`\@=11
\font\tenln    = line10
\font\tenlnw   = linew10

\newskip\Einheit \Einheit=0.5cm
\newcount\xcoord \newcount\ycoord
\newdimen\xdim \newdimen\ydim \newdimen\PfadD@cke \newdimen\Pfadd@cke

\newcount\@tempcnta
\newcount\@tempcntb

\newdimen\@tempdima
\newdimen\@tempdimb

\newdimen\@wholewidth
\newdimen\@halfwidth

\newcount\@xarg
\newcount\@yarg
\newcount\@yyarg
\newbox\@linechar
\newbox\@tempboxa
\newdimen\@linelen
\newdimen\@clnwd
\newdimen\@clnht

\newif\if@negarg

\def\@whilenoop#1{}
\def\@whiledim#1\do #2{\ifdim #1\relax#2\@iwhiledim{#1\relax#2}\fi}
\def\@iwhiledim#1{\ifdim #1\let\@nextwhile=\@iwhiledim
        \else\let\@nextwhile=\@whilenoop\fi\@nextwhile{#1}}

\def\@whileswnoop#1\fi{}
\def\@whilesw#1\fi#2{#1#2\@iwhilesw{#1#2}\fi\fi}
\def\@iwhilesw#1\fi{#1\let\@nextwhile=\@iwhilesw
         \else\let\@nextwhile=\@whileswnoop\fi\@nextwhile{#1}\fi}

\def\thinlines{\let\@linefnt\tenln \let\@circlefnt\tencirc
  \@wholewidth\fontdimen8\tenln \@halfwidth .5\@wholewidth}
\def\thicklines{\let\@linefnt\tenlnw \let\@circlefnt\tencircw
  \@wholewidth\fontdimen8\tenlnw \@halfwidth .5\@wholewidth}
\thinlines

\PfadD@cke1pt \Pfadd@cke0.5pt
\def\PfadDicke#1{\PfadD@cke#1 \divide\PfadD@cke by2 \Pfadd@cke\PfadD@cke \multiply\PfadD@cke by2}
\long\def\LOOP#1\REPEAT{\def\BODY{#1}\ITERATE}
\def\ITERATE{\BODY \let\next\ITERATE \else\let\next\relax\fi \next}
\let\REPEAT=\fi
\def\Punkt{\hbox{\raise-2pt\hbox to0pt{\hss$\ssize\bullet$\hss}}}
\def\DuennPunkt(#1,#2){\unskip
  \raise#2 \Einheit\hbox to0pt{\hskip#1 \Einheit
          \raise-2.5pt\hbox to0pt{\hss$\bullet$\hss}\hss}}
\def\NormalPunkt(#1,#2){\unskip
  \raise#2 \Einheit\hbox to0pt{\hskip#1 \Einheit
          \raise-3pt\hbox to0pt{\hss\twelvepoint$\bullet$\hss}\hss}}
\def\DickPunkt(#1,#2){\unskip
  \raise#2 \Einheit\hbox to0pt{\hskip#1 \Einheit
          \raise-4pt\hbox to0pt{\hss\fourteenpoint$\bullet$\hss}\hss}}
\def\Kreis(#1,#2){\unskip
  \raise#2 \Einheit\hbox to0pt{\hskip#1 \Einheit
          \raise-4pt\hbox to0pt{\hss\fourteenpoint$\circ$\hss}\hss}}

\def\Line@(#1,#2)#3{\@xarg #1\relax \@yarg #2\relax
\@linelen=#3\Einheit
\ifnum\@xarg =0 \@vline
  \else \ifnum\@yarg =0 \@hline \else \@sline\fi
\fi}

\def\@sline{\ifnum\@xarg< 0 \@negargtrue \@xarg -\@xarg \@yyarg -\@yarg
  \else \@negargfalse \@yyarg \@yarg \fi
\ifnum \@yyarg >0 \@tempcnta\@yyarg \else \@tempcnta -\@yyarg \fi
\ifnum\@tempcnta>6 \@badlinearg\@tempcnta0 \fi
\ifnum\@xarg>6 \@badlinearg\@xarg 1 \fi
\setbox\@linechar\hbox{\@linefnt\@getlinechar(\@xarg,\@yyarg)}%
\ifnum \@yarg >0 \let\@upordown\raise \@clnht\z@
   \else\let\@upordown\lower \@clnht \ht\@linechar\fi
\@clnwd=\wd\@linechar
\if@negarg \hskip -\wd\@linechar \def\@tempa{\hskip -2\wd\@linechar}\else
     \let\@tempa\relax \fi
\@whiledim \@clnwd <\@linelen \do
  {\@upordown\@clnht\copy\@linechar
   \@tempa
   \advance\@clnht \ht\@linechar
   \advance\@clnwd \wd\@linechar}%
\advance\@clnht -\ht\@linechar
\advance\@clnwd -\wd\@linechar
\@tempdima\@linelen\advance\@tempdima -\@clnwd
\@tempdimb\@tempdima\advance\@tempdimb -\wd\@linechar
\if@negarg \hskip -\@tempdimb \else \hskip \@tempdimb \fi
\multiply\@tempdima \@m
\@tempcnta \@tempdima \@tempdima \wd\@linechar \divide\@tempcnta \@tempdima
\@tempdima \ht\@linechar \multiply\@tempdima \@tempcnta
\divide\@tempdima \@m
\advance\@clnht \@tempdima
\ifdim \@linelen <\wd\@linechar
   \hskip \wd\@linechar
  \else\@upordown\@clnht\copy\@linechar\fi}

\def\@hline{\ifnum \@xarg <0 \hskip -\@linelen \fi
\vrule height\Pfadd@cke width \@linelen depth\Pfadd@cke
\ifnum \@xarg <0 \hskip -\@linelen \fi}

\def\@getlinechar(#1,#2){\@tempcnta#1\relax\multiply\@tempcnta 8
\advance\@tempcnta -9 \ifnum #2>0 \advance\@tempcnta #2\relax\else
\advance\@tempcnta -#2\relax\advance\@tempcnta 64 \fi
\char\@tempcnta}

\def\Vektor(#1,#2)#3(#4,#5){\unskip\leavevmode
  \xcoord#4\relax \ycoord#5\relax
      \raise\ycoord \Einheit\hbox to0pt{\hskip\xcoord \Einheit
         \Vector@(#1,#2){#3}\hss}}

\def\Vector@(#1,#2)#3{\@xarg #1\relax \@yarg #2\relax
\@tempcnta \ifnum\@xarg<0 -\@xarg\else\@xarg\fi
\ifnum\@tempcnta<5\relax
\@linelen=#3\Einheit
\ifnum\@xarg =0 \@vvector
  \else \ifnum\@yarg =0 \@hvector \else \@svector\fi
\fi
\else\@badlinearg\fi}

\def\@hvector{\@hline\hbox to 0pt{\@linefnt
\ifnum \@xarg <0 \@getlarrow(1,0)\hss\else
    \hss\@getrarrow(1,0)\fi}}

\def\@vvector{\ifnum \@yarg <0 \@downvector \else \@upvector \fi}

\def\@svector{\@sline
\@tempcnta\@yarg \ifnum\@tempcnta <0 \@tempcnta=-\@tempcnta\fi
\ifnum\@tempcnta <5
  \hskip -\wd\@linechar
  \@upordown\@clnht \hbox{\@linefnt  \if@negarg
  \@getlarrow(\@xarg,\@yyarg) \else \@getrarrow(\@xarg,\@yyarg) \fi}%
\else\@badlinearg\fi}

\def\@upline{\hbox to \z@{\hskip -.5\Pfadd@cke \vrule width \Pfadd@cke
   height \@linelen depth \z@\hss}}

\def\@downline{\hbox to \z@{\hskip -.5\Pfadd@cke \vrule width \Pfadd@cke
   height \z@ depth \@linelen \hss}}

\def\@upvector{\@upline\setbox\@tempboxa\hbox{\@linefnt\char'66}\raise
     \@linelen \hbox to\z@{\lower \ht\@tempboxa\box\@tempboxa\hss}}

\def\@downvector{\@downline\lower \@linelen
      \hbox to \z@{\@linefnt\char'77\hss}}

\def\@getlarrow(#1,#2){\ifnum #2 =\z@ \@tempcnta='33\else
\@tempcnta=#1\relax\multiply\@tempcnta \sixt@@n \advance\@tempcnta
-9 \@tempcntb=#2\relax\multiply\@tempcntb \tw@
\ifnum \@tempcntb >0 \advance\@tempcnta \@tempcntb\relax
\else\advance\@tempcnta -\@tempcntb\advance\@tempcnta 64
\fi\fi\char\@tempcnta}

\def\@getrarrow(#1,#2){\@tempcntb=#2\relax
\ifnum\@tempcntb < 0 \@tempcntb=-\@tempcntb\relax\fi
\ifcase \@tempcntb\relax \@tempcnta='55 \or
\ifnum #1<3 \@tempcnta=#1\relax\multiply\@tempcnta
24 \advance\@tempcnta -6 \else \ifnum #1=3 \@tempcnta=49
\else\@tempcnta=58 \fi\fi\or
\ifnum #1<3 \@tempcnta=#1\relax\multiply\@tempcnta
24 \advance\@tempcnta -3 \else \@tempcnta=51\fi\or
\@tempcnta=#1\relax\multiply\@tempcnta
\sixt@@n \advance\@tempcnta -\tw@ \else
\@tempcnta=#1\relax\multiply\@tempcnta
\sixt@@n \advance\@tempcnta 7 \fi\ifnum #2<0 \advance\@tempcnta 64 \fi
\char\@tempcnta}

\def\Diagonale(#1,#2)#3{\unskip\leavevmode
  \xcoord#1\relax \ycoord#2\relax
      \raise\ycoord \Einheit\hbox to0pt{\hskip\xcoord \Einheit
         \Line@(1,1){#3}\hss}}
\def\AntiDiagonale(#1,#2)#3{\unskip\leavevmode
  \xcoord#1\relax \ycoord#2\relax 
      \raise\ycoord \Einheit\hbox to0pt{\hskip\xcoord \Einheit
         \Line@(1,-1){#3}\hss}}
\def\Pfad(#1,#2),#3\endPfad{\unskip\leavevmode
  \xcoord#1 \ycoord#2 \thicklines\ZeichnePfad#3\endPfad\thinlines}
\def\ZeichnePfad#1{\ifx#1\endPfad\let\next\relax
  \else\let\next\ZeichnePfad
    \ifnum#1=1
      \raise\ycoord \Einheit\hbox to0pt{\hskip\xcoord \Einheit
         \vrule height\Pfadd@cke width1 \Einheit depth\Pfadd@cke\hss}%
      \advance\xcoord by 1
    \else\ifnum#1=2
      \raise\ycoord \Einheit\hbox to0pt{\hskip\xcoord \Einheit
        \hbox{\hskip-\PfadD@cke\vrule height1 \Einheit width\PfadD@cke depth0pt}\hss}%
      \advance\ycoord by 1
    \else\ifnum#1=3
      \raise\ycoord \Einheit\hbox to0pt{\hskip\xcoord \Einheit
         \Line@(1,1){1}\hss}
      \advance\xcoord by 1
      \advance\ycoord by 1
    \else\ifnum#1=4
      \raise\ycoord \Einheit\hbox to0pt{\hskip\xcoord \Einheit
         \Line@(1,-1){1}\hss}
      \advance\xcoord by 1
      \advance\ycoord by -1
    \else\ifnum#1=5
      \advance\xcoord by -1
      \raise\ycoord \Einheit\hbox to0pt{\hskip\xcoord \Einheit
         \vrule height\Pfadd@cke width1 \Einheit depth\Pfadd@cke\hss}%
    \else\ifnum#1=6
      \advance\ycoord by -1
      \raise\ycoord \Einheit\hbox to0pt{\hskip\xcoord \Einheit
        \hbox{\hskip-\PfadD@cke\vrule height1 \Einheit width\PfadD@cke depth0pt}\hss}%
    \else\ifnum#1=7
      \advance\xcoord by -1
      \advance\ycoord by -1
      \raise\ycoord \Einheit\hbox to0pt{\hskip\xcoord \Einheit
         \Line@(1,1){1}\hss}
    \else\ifnum#1=8
      \advance\xcoord by -1
      \advance\ycoord by +1
      \raise\ycoord \Einheit\hbox to0pt{\hskip\xcoord \Einheit
         \Line@(1,-1){1}\hss}
    \fi\fi\fi\fi
    \fi\fi\fi\fi
  \fi\next}
\def\hSSchritt{\leavevmode\raise-.4pt\hbox to0pt{\hss.\hss}\hskip.2\Einheit
  \raise-.4pt\hbox to0pt{\hss.\hss}\hskip.2\Einheit
  \raise-.4pt\hbox to0pt{\hss.\hss}\hskip.2\Einheit
  \raise-.4pt\hbox to0pt{\hss.\hss}\hskip.2\Einheit
  \raise-.4pt\hbox to0pt{\hss.\hss}\hskip.2\Einheit}
\def\vSSchritt{\vbox{\baselineskip.2\Einheit\lineskiplimit0pt
\hbox{.}\hbox{.}\hbox{.}\hbox{.}\hbox{.}}}
\def\DSSchritt{\leavevmode\raise-.4pt\hbox to0pt{%
  \hbox to0pt{\hss.\hss}\hskip.2\Einheit
  \raise.2\Einheit\hbox to0pt{\hss.\hss}\hskip.2\Einheit
  \raise.4\Einheit\hbox to0pt{\hss.\hss}\hskip.2\Einheit
  \raise.6\Einheit\hbox to0pt{\hss.\hss}\hskip.2\Einheit
  \raise.8\Einheit\hbox to0pt{\hss.\hss}\hss}}
\def\dSSchritt{\leavevmode\raise-.4pt\hbox to0pt{%
  \hbox to0pt{\hss.\hss}\hskip.2\Einheit
  \raise-.2\Einheit\hbox to0pt{\hss.\hss}\hskip.2\Einheit
  \raise-.4\Einheit\hbox to0pt{\hss.\hss}\hskip.2\Einheit
  \raise-.6\Einheit\hbox to0pt{\hss.\hss}\hskip.2\Einheit
  \raise-.8\Einheit\hbox to0pt{\hss.\hss}\hss}}
\def\SPfad(#1,#2),#3\endSPfad{\unskip\leavevmode
  \xcoord#1 \ycoord#2 \ZeichneSPfad#3\endSPfad}
\def\ZeichneSPfad#1{\ifx#1\endSPfad\let\next\relax
  \else\let\next\ZeichneSPfad
    \ifnum#1=1
      \raise\ycoord \Einheit\hbox to0pt{\hskip\xcoord \Einheit
         \hSSchritt\hss}%
      \advance\xcoord by 1
    \else\ifnum#1=2
      \raise\ycoord \Einheit\hbox to0pt{\hskip\xcoord \Einheit
        \hbox{\hskip-2pt \vSSchritt}\hss}%
      \advance\ycoord by 1
    \else\ifnum#1=3
      \raise\ycoord \Einheit\hbox to0pt{\hskip\xcoord \Einheit
         \DSSchritt\hss}
      \advance\xcoord by 1
      \advance\ycoord by 1
    \else\ifnum#1=4
      \raise\ycoord \Einheit\hbox to0pt{\hskip\xcoord \Einheit
         \dSSchritt\hss}
      \advance\xcoord by 1
      \advance\ycoord by -1
    \else\ifnum#1=5
      \advance\xcoord by -1
      \raise\ycoord \Einheit\hbox to0pt{\hskip\xcoord \Einheit
         \hSSchritt\hss}%
    \else\ifnum#1=6
      \advance\ycoord by -1
      \raise\ycoord \Einheit\hbox to0pt{\hskip\xcoord \Einheit
        \hbox{\hskip-2pt \vSSchritt}\hss}%
    \else\ifnum#1=7
      \advance\xcoord by -1
      \advance\ycoord by -1
      \raise\ycoord \Einheit\hbox to0pt{\hskip\xcoord \Einheit
         \DSSchritt\hss}
    \else\ifnum#1=8
      \advance\xcoord by -1
      \advance\ycoord by 1
      \raise\ycoord \Einheit\hbox to0pt{\hskip\xcoord \Einheit
         \dSSchritt\hss}
    \fi\fi\fi\fi
    \fi\fi\fi\fi
  \fi\next}
\def\Koordinatenachsen(#1,#2){\unskip
 \hbox to0pt{\hskip-.5pt\vrule height#2 \Einheit width.5pt depth1 \Einheit}%
 \hbox to0pt{\hskip-1 \Einheit \xcoord#1 \advance\xcoord by1
    \vrule height0.25pt width\xcoord \Einheit depth0.25pt\hss}}
\def\Koordinatenachsen(#1,#2)(#3,#4){\unskip
 \hbox to0pt{\hskip-.5pt \ycoord-#4 \advance\ycoord by1
    \vrule height#2 \Einheit width.5pt depth\ycoord \Einheit}%
 \hbox to0pt{\hskip-1 \Einheit \hskip#3\Einheit 
    \xcoord#1 \advance\xcoord by1 \advance\xcoord by-#3 
    \vrule height0.25pt width\xcoord \Einheit depth0.25pt\hss}}
\def\Gitter(#1,#2){\unskip \xcoord0 \ycoord0 \leavevmode
  \LOOP\ifnum\ycoord<#2
    \loop\ifnum\xcoord<#1
      \raise\ycoord \Einheit\hbox to0pt{\hskip\xcoord \Einheit\Punkt\hss}%
      \advance\xcoord by1
    \repeat
    \xcoord0
    \advance\ycoord by1
  \REPEAT}
\def\Gitter(#1,#2)(#3,#4){\unskip \xcoord#3 \ycoord#4 \leavevmode
  \LOOP\ifnum\ycoord<#2
    \loop\ifnum\xcoord<#1
      \raise\ycoord \Einheit\hbox to0pt{\hskip\xcoord \Einheit\Punkt\hss}%
      \advance\xcoord by1
    \repeat
    \xcoord#3
    \advance\ycoord by1
  \REPEAT}
\def\Label#1#2(#3,#4){\unskip \xdim#3 \Einheit \ydim#4 \Einheit
  \def\lo{\advance\xdim by-.5 \Einheit \advance\ydim by.5 \Einheit}%
  \def\llo{\advance\xdim by-.25cm \advance\ydim by.5 \Einheit}%
  \def\loo{\advance\xdim by-.5 \Einheit \advance\ydim by.25cm}%
  \def\o{\advance\ydim by.25cm}%
  \def\ro{\advance\xdim by.5 \Einheit \advance\ydim by.5 \Einheit}%
  \def\rro{\advance\xdim by.25cm \advance\ydim by.5 \Einheit}%
  \def\roo{\advance\xdim by.5 \Einheit \advance\ydim by.25cm}%
  \def\l{\advance\xdim by-.30cm}%
  \def\r{\advance\xdim by.30cm}%
  \def\lu{\advance\xdim by-.5 \Einheit \advance\ydim by-.6 \Einheit}%
  \def\llu{\advance\xdim by-.25cm \advance\ydim by-.6 \Einheit}%
  \def\luu{\advance\xdim by-.5 \Einheit \advance\ydim by-.30cm}%
  \def\u{\advance\ydim by-.30cm}%
  \def\ru{\advance\xdim by.5 \Einheit \advance\ydim by-.6 \Einheit}%
  \def\rru{\advance\xdim by.25cm \advance\ydim by-.6 \Einheit}%
  \def\ruu{\advance\xdim by.5 \Einheit \advance\ydim by-.30cm}%
  #1\raise\ydim\hbox to0pt{\hskip\xdim
     \vbox to0pt{\vss\hbox to0pt{\hss$#2$\hss}\vss}\hss}%
}
\catcode`\@=13

\def\BereAA{1}
\def\BrFoAA{2}
\def\BurrAA{3}
\def\BuMMAA{4}
\def\BuCFAA{5}
\def\ChDDAB{6}
\def\DeDFAA{7}
\def\DuSaAA{8}
\def\FomiAZ{9}
\def\FomiAB{10}
\def\FomiAF{11}
\def\GeZeAA{12}
\def\GoulAD{13}
\def\GrMrAA{14}
\def\GreCAA{15}
\def\KnutAA{16}
\def\KratBC{17}
\def\KratCE{18}
\def\LeeuAD{19}
\def\LeeuAH{20}
\def\MacdAC{21}
\def\ProcAA{22}
\def\RobiAA{23}
\def\RobyAA{24}
\def\RobyAD{25}
\def\SagaAQ{26}
\def\ScheAA{27}
\def\SchuAB{28}
\def\SchuAA{29}
\def\StanBI{30}
\def\SunaAE{31}
\def\SunaAD{32}
\def\SunaAC{33}

\def\SB{2}
\def\SC{3}

\def\TA{1}
\def\TB{2}
\def\TC{3}
\def\TD{4}

\def\FA{1}
\def\FB{2}
\def\FC{3}
\def\FD{4}
\def\FE{5}
\def\FF{6}
\def\FG{7}

\def\AA{4.1}

\topmatter 
\title Bijections between oscillating tableaux and (semi)standard
tableaux via growth diagrams
\endtitle 
\author C.~Krattenthaler 
\endauthor 
\affil 
Fakult\"at f\"ur Mathematik, Universit\"at Wien,\\
Oskar-Morgenstern-Platz~1, A-1090 Vienna, Austria.\\
WWW: {\tt http://www.mat.univie.ac.at/\~{}kratt}
\endaffil
\address Fakult\"at f\"ur Mathematik, Universit\"at Wien,
Oskar-Morgenstern-Platz~1, A-1090 Vienna,\linebreak Austria.\newline
WWW: \tt http://www.mat.univie.ac.at/\~{}kratt
\endaddress

\thanks Research partially supported 
by the Austrian Science Foundation FWF, grant Z130-N13, and
grant SFB F50 (Special Research Program
``Algorithmic and Enumerative Combinatorics")
\endthanks
\subjclass Primary 05A15;
 Secondary 05E10
\endsubjclass
\keywords standard Young tableaux, semistandard tableaux,
oscillating tableaux,
Robin\-son--Schensted correspondence, 
Robin\-son--Schensted--Knuth correspondence, 
growth diagrams
\endkeywords
\abstract 
We prove that 
the number of oscillating tableaux
of length $n$ with at most $k$ columns,
starting at $\emptyset$ and ending at the one-column shape $(1^m)$,
is equal to the number of standard Young tableaux of size~$n$
with $m$ columns of odd length, all columns of length
at most $2k$. This refines a conjecture of Burrill, 
which it thereby establishes. 
We prove as well a ``Knuth-type" extension stating a similar
equi-enumeration result between generalised oscillating tableaux
and semistandard tableaux.
\endabstract
\endtopmatter
\document

\subhead 1. Introduction \endsubhead
The {\it Robinson--Schensted correspondence} \cite{\RobiAA, \ScheAA}
(see \cite{\SagaAQ, Sec.~3.1}) is a bijection between {\it permutations}
of $\{1,2,\dots,n\}$ and pairs of {\it standard (Young) tableaux} of the
same shape of size $n$ (see Section~\SB\ for all definitions).
Knuth's extension, the so-called 
{\it Robinson--Schensted--Knuth} ({\it RSK}) {\it correspondence} 
\cite{\KnutAA}
(see \cite{\SagaAQ, Sec.~4.8}) is a bijection between {\it
non-negative integer matrices}
and pairs of {\it semistandard tableaux} of the
same shape. These correspondences are not only attractive in
their own right due to
their elegance, but are important since they map several natural
statistics of permutations, respectively matrices, to corresponding ones for
tableaux, and thus allow for numerous refinements. 
An important statistic in this context is the length of the
longest increasing (or decreasing) subsequence in a permutation
(or of certain chains of entries in a matrix) which, by Schensted's
theorem \cite{\ScheAA}, are mapped to the length of the first row
(or first column) of the shapes of the tableaux. Greene \cite{\GreCAA}
extended Schensted's theorem by describing precisely how lengths 
of increasing (or decreasing) subsequences in permutations (chains in
matrices) determine the shape of the tableaux in the image pair under
these correspondences.

Standard or semistandard 
tableaux may be seen as sequences of Ferrers diagrams, where 
an element in the sequence is followed by a Ferrers diagram which
is by one cell (in the case of standard tableaux) or by a horizontal
strip (in the case of semistandard tableaux) larger. A variation
consists in considering sequences of Ferrers diagrams where, from
one element in the sequence to the next, one also allows sometimes
to shrink by a cell or by a horizontal (or vertical) strip. This leads to the
notion of {\it oscillating tableaux}. Also in this context, there
are Robinson--Schensted(--Knuth) like algorithms which connect
oscillating tableaux to, for instance, 
involutions, matchings, or set partitions;
see \cite{\BereAA, \DeDFAA, \DuSaAA, 
\ProcAA, \RobyAA, \RobyAD, \SunaAE, \SunaAD, \SunaAC}.
Greene's theorem \cite{\GreCAA} still applies, which has been
particularly exploited in \cite{\ChDDAB, \KratCE}.

If one encounters families of tableaux, permutations, integer
matrices, etc. which seem to be enumerated by the same numbers,
then one must suspect that an RSK-like bijection lurks in the
background. The purpose of the present paper is 
to illustrate this ``principle" by
applying it to a recent conjecture of Burrill
\cite{\BurrAA, Conj.~6.2.1} (see 
\cite{\BuMMAA, Conj.~4} for the formulation below; 
again, for all undefined terminology
see Section~\SB).

\proclaim{Conjecture \smc(Burrill)}
Let $n$ and $k$ be given non-negative integers. 
The number of oscillating tableaux
of length $n$ with at most $k$ columns,
starting at $\emptyset$ and ending at some one-column shape
is the same as the number of
standard Young tableaux of size~$n$ (meaning that its shape has $n$~cells) 
with all columns of length at most $2k$.
\endproclaim

In Theorem~\TC\ in Section~\SC, a refinement of this conjecture will
be established, which in addition relates the length of the one-column shape
in which the oscillating tableaux end to the number of odd-length
columns of the standard tableaux. Moreover, in Theorem~\TD,
we present a ``Knuth-type" extension, in which we allow
more general oscillating tableaux, and the standard Young tableaux
get replaced by semistandard tableaux. We point out that a
different bijection proving Theorem~\TC\ has been presented
by Burrill, Courtiel, Fusy, Melczer and Mishna in \cite{\BuCFAA},
the difference lying in the way the parameter~$m$ (the number of
odd columns of the standard Young tableau, respectively the
size of the final shape of the oscillating tableau) is kept track of.

While, originally, RSK-like correspondences are based on
insertion-deletion algorithms, it is nowadays standard that the
most transparent way to present these correspondences is by
means of Fomin's {\it growth
diagrams} \cite{\FomiAZ, \FomiAB, \FomiAF}
(see \cite{\RobyAA, \RobyAD},
\cite{\SagaAQ, Sec.~5.2} and \cite{\StanBI, Sec.~7.13} for non-technical
expositions). This is also the point of view we shall adopt in
our (bijective) proofs of Theorems~\TC\ and \TD. It will be combined
with an application of Sch\"utzenberger's \cite{\SchuAA} jeu de
taquin (for which also geometric realisations have been proposed ---
see \cite{\LeeuAD} ---, which we shall however not use here).

In the final Section~4, we explain that --- non-illuminating ---
computational proofs of Theorems~\TC\ and \TD\ can be extracted
from the literature by appropriately combining results of Gessel and
Zeilberger \cite{\GeZeAA} and of Goulden \cite{\GoulAD}, we
comment on what happens if, instead of an even number, we bound
the length of columns of the standard Young tableaux in Theorem~\TC\
by an odd number, and we provide a more detailed discussion of
the differences between the bijection proving Theorem~\TC\ presented
here and the one in \cite{\BuCFAA}.

\subhead 2. Definitions and notation \endsubhead
We start by fixing the standard partition notation (cf.\
e.g.\ \cite{\StanBI, Sec.~7.2}).
A {\it partition} is a weakly decreasing sequence
$\la=(\la_1,\la_2,\dots,\la_\ell)$ of positive integers.
This also includes the {\it empty partition} $()$, denoted by
$\emptyset$. For the sake of convenience, we shall often tacitly
identify a partition $\la=(\la_1,\la_2,\dots,\la_\ell)$ with the
infinite sequence $(\la_1,\la_2,\dots,\la_\ell,0,0,\dots)$, that is,
the sequence which arises from $\la$ by appending infinitely
many $0$'s. To each partition $\la$, one
associates its {\it Ferrers diagram} (also called {\it Ferrers shape}), 
which is the left-justified
arrangement of squares with $\la_i$ squares in the $i$-th row,
$i=1,2,\dots$. The number of squares in the Ferrers diagram,
$\la_1+\la_2+\dots+\la_\ell$, is called the {\it size} of the
partition~$\la$, and is denoted by $\vert\la\vert$.
We define a {\it partial order} $\subseteq$ on partitions 
by containment of their Ferrers diagrams.
The {\it union} $\mu\cup\nu$ of two partitions $\mu$ and $\nu$
is the partition which arises by forming the union of the Ferrers
diagrams of $\mu$ and $\nu$. Thus, if
$\mu=(\mu_1,\mu_2,\dots)$ and
$\nu=(\nu_1,\nu_2,\dots)$, then $\mu\cup\nu$ is the partition 
$\la=(\la_1,\la_2,\dots)$, where $\la_i=\max\{\mu_i,\nu_i\}$ for 
$i=1,2,\dots$. The {\it intersection} 
$\mu\cap\nu$ of two partitions $\mu$ and $\nu$
is the partition which arises by forming the intersection of the Ferrers
diagrams of $\mu$ and $\nu$. Thus, if
$\mu=(\mu_1,\mu_2,\dots)$ and
$\nu=(\nu_1,\nu_2,\dots)$, then $\mu\cap\nu$ is the partition 
$\rh=(\rh_1,\rh_2,\dots)$, where $\rh_i=\min\{\mu_i,\nu_i\}$ for 
$i=1,2,\dots$.

\medskip
Given a partition $\la=(\la_1,\la_2,\dots,\la_\ell)$, a {\it standard
(Young) tableau of shape\/} $\la$ is a left-justified arrangement
of positive integers with $\la_i$ entries in row~$i$, $i=1,2,\dots$,
such that the entries along rows and columns are increasing.
An arrangement $T$ of the same form as above is called
{\it semistandard tableau of shape\/} $\la$ if 
the entries along rows are weakly increasing
and such that the entries along columns are strictly increasing. 
By considering the sequence of partitions (Ferrers shapes)
$(\la^{i})_{i\ge0}$, where $\la^i$ is the shape formed by the
entries of $T$ which are at most~$i$,
$i=0,1,2,\dots$, one sees that standard tableaux of shape $\la$ are in
bijection with sequences 
$\emptyset=\la^0\subseteq\la^1\subseteq\dots\subseteq\la^n=\la$,
where $\la^{i-1}$ and $\la^i$ differ by exactly one square for all~$i$,
while semistandard tableaux of shape $\la$ are in bijection with
such sequences where $\la^{i-1}$ and $\la^i$ differ by a 
{\it horizontal strip} for all~$i$, that is, by an arrangement of squares with
at most one square in each column.

Generalising the above concepts, we call a sequence
$\emptyset=\la^0,\la^1,\dots,\la^n=\la$ of partitions an
{\it oscillating tableau of shape} $\la$ if either $\la^{i-1}\subseteq
\la^i$ or $\la^{i-1}\supseteq \la^i$ and $\la^{i-1}$ and $\la^i$
differ by exactly one square, $i=1,2,\dots,n$. 
The number~$n$ is called the {\it length\/} of the (generalised)
oscillating tableau.
If we say that an oscillating tableau has ``at most $k$ 
columns" then we mean that {\it all\/} partitions in the sequence
have at most $k$ columns.

\medskip
{\it Growth diagrams} are certain labellings of arrangements of cells.
The arrangements of cells which we need here are arrangements which
are left-justified (that is, they have a straight vertical left
boundary),
bottom-justified (that is, they have a straight horizontal bottom
boundary), and rows and columns in the arrangement are ``without"
holes, that is, if we move along the top-right boundary of the
arrangement, we always move either to the right or to the bottom.
Figure~\FA.a shows an example of such a cell arrangement.

\midinsert
$$
\Einheit.4cm
\Pfad(0,18),666666666666666666111111111111111\endPfad
\Pfad(0,18),111666666666666666666\endPfad
\Pfad(0,15),111111666666666666666\endPfad
\Pfad(0,12),111111666111666666666\endPfad
\Pfad(0,9),111111111666666111666\endPfad
\Pfad(0,6),111111111666111111666\endPfad
\Pfad(0,3),111111111111111666\endPfad
\hbox{\hskip7.5cm}
\Pfad(0,18),666666666666666666111111111111111\endPfad
\Pfad(0,18),111666666666666666666\endPfad
\Pfad(0,15),111111666666666666666\endPfad
\Pfad(0,12),111111666111666666666\endPfad
\Pfad(0,9),111111111666666111666\endPfad
\Pfad(0,6),111111111666111111666\endPfad
\Pfad(0,3),111111111111111666\endPfad
\Label\ro{\text {\seventeenpoint $0$}}(1,1)
\Label\ro{\text {\seventeenpoint $0$}}(4,1)
\Label\ro{\text {\seventeenpoint $0$}}(7,1)
\Label\ro{\text {\seventeenpoint $0$}}(10,1)
\Label\ro{\text {\seventeenpoint $0$}}(1,4)
\Label\ro{\text {\seventeenpoint $0$}}(7,4)
\Label\ro{\text {\seventeenpoint $0$}}(1,7)
\Label\ro{\text {\seventeenpoint $0$}}(4,7)
\Label\ro{\text {\seventeenpoint $0$}}(7,7)
\Label\ro{\text {\seventeenpoint $0$}}(4,10)
\Label\ro{\text {\seventeenpoint $0$}}(1,13)
\Label\ro{\text {\seventeenpoint $0$}}(4,13)
\Label\ro{\text {\seventeenpoint $0$}}(1,16)
\Label\ro{\text {\seventeenpoint $1$}}(4,4)
\Label\ro{\text {\seventeenpoint $1$}}(1,10)
\Label\ro{\text {\seventeenpoint $1$}}(13,1)
\hskip6cm
$$
\vskip10pt
\centerline{\eightpoint a. A cell arrangement
\hskip4cm
b. A filling of the cell arrangement}
\vskip6pt
\centerline{\smc Figure \FA}
\endinsert

We fill the cells of such an arrangement $C$ with non-negative integers.
Most of the time,
the fillings will be restricted to $0$-$1$-fillings
such that every row and every column contains at most one $1$.
See Figure~\FA.b for an example.

Next, the corners of the cells are labelled by partitions
such that the following two conditions are satisfied:

\roster
\item"(C1)" A partition is either equal to its right neighbour
or smaller by exactly one square, the same being true for a
partition and its top neighbour.
\item"(C2)" A partition and its right neighbour are equal if and only
if in the column of cells of $C$ below them there appears no $1$ and if
their bottom neighbours are also equal to each other.
Similarly, a partition and its top neighbour are equal if and only if
in the row of cells of $C$ to the left of them there appears no $1$ and if
their left neighbours are also equal to each other.
\endroster

See Figure~\FB\ for an example. (More examples can be found in
Figures~\FD--\FF.)
There, we use a short notation for partitions. For example,
$11$ is short for $(1,1)$.
Moreover, we changed the convention of representing the filling slightly
for better visibility,
by suppressing $0$'s and by replacing $1$'s by X's. Indeed, the
filling represented in Figure~\FB\ is the same as the one in
Figure~\FA.b.

Diagrams which obey the conditions (C1) and (C2) are called {\it
growth diagrams}.

\vskip10pt
\midinsert
$$
\Einheit.4cm
\Pfad(0,18),666666666666666666111111111111111\endPfad
\Pfad(0,18),111666666666666666666\endPfad
\Pfad(0,15),111111666666666666666\endPfad
\Pfad(0,12),111111666111666666666\endPfad
\Pfad(0,9),111111111666666111666\endPfad
\Pfad(0,6),111111111666111111666\endPfad
\Pfad(0,3),111111111111111666\endPfad
\Label\ro{\emptyset}(0,0)
\Label\ro{\emptyset}(0,3)
\Label\ro{\emptyset}(0,6)
\Label\ro{\emptyset}(0,9)
\Label\ro{\emptyset}(0,12)
\Label\ro{\emptyset}(0,15)
\Label\ro{\emptyset}(0,18)
\Label\ro{1}(3,18)
\Label\ro{1}(3,15)
\Label\ro{1}(3,12)
\Label\ro{\emptyset}(3,9)
\Label\ro{\emptyset}(3,6)
\Label\ro{\emptyset}(3,3)
\Label\ro{\emptyset}(3,0)
\Label\ro{\hphantom{1}11}(6,15)
\Label\ro{\hphantom{1}11}(6,12)
\Label\ro{1}(6,9)
\Label\ro{1}(6,6)
\Label\ro{\emptyset}(6,3)
\Label\ro{\emptyset}(6,0)
\Label\ro{1}(9,9)
\Label\ro{1}(9,6)
\Label\ro{\emptyset}(9,3)
\Label\ro{\emptyset}(9,0)
\Label\ro{\emptyset}(12,3)
\Label\ro{\emptyset}(12,0)
\Label\ro{1}(15,3)
\Label\ro{\emptyset}(15,0)
\Label\ro{\text {\seventeenpoint X}}(4,4)
\Label\ro{\text {\seventeenpoint X}}(1,10)
\Label\ro{\text {\seventeenpoint X}}(13,1)
\hskip6cm
$$
\vskip10pt
\centerline{\eightpoint A growth diagram}
\vskip6pt
\centerline{\smc Figure \FB}
\endinsert
\vskip10pt

We are interested in growth diagrams which
obey the following ({\it forward}) {\it local rules} 
(see Figure~\FC).

\vskip10pt
\midinsert
$$
\Pfad(0,0),1111222255556666\endPfad
\Label\lu{\rh}(0,0)
\Label\ru{\mu}(4,0)
\Label\lo{\nu}(0,4)
\Label\ro{\la}(4,4)
\hbox{\hskip6cm}
\Pfad(0,0),1111222255556666\endPfad
\Label\lu{\rh}(0,0)
\Label\ru{\mu}(4,0)
\Label\lo{\nu}(0,4)
\Label\ro{\la}(4,4)
\thicklines
\Pfad(1,1),33\endPfad
\raise.5pt\hbox to 0pt{\hbox{\hskip-.5pt}\Pfad(1,1),3\endPfad\hss}%
\raise.5pt\hbox to 0pt{\hbox{\hskip-.5pt}\Pfad(2,2),3\endPfad\hss}%
\raise1pt\hbox to 0pt{\hbox{\hskip-1pt}\Pfad(1,1),3\endPfad\hss}%
\raise1pt\hbox to 0pt{\hbox{\hskip-1pt}\Pfad(2,2),3\endPfad\hss}%
\raise-.5pt\hbox to 0pt{\hbox{\hskip.5pt}\Pfad(1,1),3\endPfad\hss}%
\raise-.5pt\hbox to 0pt{\hbox{\hskip.5pt}\Pfad(2,2),3\endPfad\hss}%
\raise-1pt\hbox to 0pt{\hbox{\hskip1pt}\Pfad(1,1),3\endPfad\hss}%
\raise-1pt\hbox to 0pt{\hbox{\hskip1pt}\Pfad(2,2),3\endPfad\hss}%
\Pfad(1,3),44\endPfad
\raise.5pt\hbox to 0pt{\hbox{\hskip.5pt}\Pfad(1,3),4\endPfad\hss}%
\raise.5pt\hbox to 0pt{\hbox{\hskip.5pt}\Pfad(2,2),4\endPfad\hss}%
\raise1pt\hbox to 0pt{\hbox{\hskip1pt}\Pfad(1,3),4\endPfad\hss}%
\raise1pt\hbox to 0pt{\hbox{\hskip1pt}\Pfad(2,2),4\endPfad\hss}%
\raise-.5pt\hbox to 0pt{\hbox{\hskip-.5pt}\Pfad(1,3),4\endPfad\hss}%
\raise-.5pt\hbox to 0pt{\hbox{\hskip-.5pt}\Pfad(2,2),4\endPfad\hss}%
\raise-1pt\hbox to 0pt{\hbox{\hskip-1pt}\Pfad(1,3),4\endPfad\hss}%
\raise-1pt\hbox to 0pt{\hbox{\hskip-1pt}\Pfad(2,2),4\endPfad\hss}%
\hskip2cm
$$
\vskip10pt
\centerline{\eightpoint a. A cell without cross\hskip3cm
b. A cell with cross}
\vskip6pt
\centerline{\smc Figure \FC}
\endinsert
\vskip10pt

\roster
\item"(F1)" If $\rh=\mu=\nu$, and if there is no cross in the
cell, then $\la=\rh$.
\item"(F2)" If $\rh=\mu\ne\nu$, then $\la=\nu$.
\item"(F3)" If $\rh=\nu\ne\mu$, then $\la=\mu$.
\item"(F4)" If $\rh,\mu,\nu$ are pairwise different, 
then $\la=\mu\cup\nu$.
\item"(F5)" If $\rh\ne \mu=\nu$, 
then $\la$ is formed by adding a square to the $(k+1)$-st row of
$\mu=\nu$, given that $\mu=\nu$ and $\rh$ differ in the $k$-th row.
\item"(F6)" If $\rh=\mu=\nu$, and if there is a cross in the
cell, then $\la$ is formed by adding a square to the first row of
$\rh=\mu=\nu$.
\endroster

Thus, if we label all the corners along the left and the bottom boundary
by empty partitions (which we shall always do in this paper), 
these rules allow one to determine all other labels
of corners uniquely.

It is not difficult to see that
the rules (F5) and (F6) are designed so that
one can also work one's way in the other direction, that is, given
$\la,\mu,\nu$, one can reconstruct $\rh$ {\it and\/} the filling of the cell.
The corresponding ({\it backward}) {\it local rules} are:

\roster
\item"(B1)" If $\la=\mu=\nu$, then $\rh=\la$.
\item"(B2)" If $\la=\mu\ne\nu$, then $\rh=\nu$.
\item"(B3)" If $\la=\nu\ne\mu$, then $\rh=\mu$.
\item"(B4)" If $\la,\mu,\nu$ are pairwise different, 
then $\rh=\mu\cap\nu$.
\item"(B5)" If $\la\ne \mu=\nu$, 
then $\rh$ is formed by deleting a square from the $(k-1)$-st row of
$\mu=\nu$, given that $\mu=\nu$ and $\la$ differ in the $k$-th row,
$k\ge2$.
\item"(B6)" If $\la\ne \mu=\nu$, and if $\la$ and $\mu=\nu$ differ in
the first row, then $\rh=\mu=\nu$.
\item""\hskip-36pt\vbox{\hsize13.76cm\noindent
In case (B6) the cell is filled with a $1$ (an X).
In all other cases the cell is filled with a $0$.}
\endroster

Thus, given a labelling of the corners along the right/up boundary of a
cell arrangement, one can algorithmically reconstruct the labels of
the other corners of the cells {\it and\/} of the $0$-$1$-filling 
by working one's way to the left and to the bottom.
These observations lead to the following theorem.

\proclaim{Theorem \TA} 
Let $C$ be an arrangement of cells.
The $0$-$1$-fillings of $C$ with the property that every row and
every column contains at most one $1$ are in bijection with labellings
$(\emptyset=\la^0,\la^1,\dots,\la^{k}=\emptyset)$ of the corners of
cells appearing along the top-right boundary of $C$, where $\la^{i-1}$
and $\la^i$ differ by at most one square, and
$\la^{i-1}\subseteq\la^i$ if $\la^{i-1}$ and $\la^i$ appear along
a horizontal edge, whereas
$\la^{i-1}\supseteq\la^i$ if $\la^{i-1}$ and $\la^i$ appear along
a vertical edge.
Moreover, $\la^{i-1}\subsetneqq\la^i$ if and only if there is a $1$
in the column of cells of $C$ below the corners labelled by $\la^{i-1}$ and
$\la^i$, and $\la^{i-1}\supsetneqq\la^i$ if and only if there is a $1$
in the row of cells of $C$ to the left of the corners labelled by
$\la^{i-1}$ and $\la^i$.
\endproclaim

\medskip
In addition to its local description,
the bijection of the above theorem has also a {\it global\/}
description. The latter is a consequence of a theorem of Greene
\cite{\GreCAA} (see also \cite{\BrFoAA, Theorems~2.1 and 3.2}). 
In order to formulate the
result, we need the following definitions: a {\it NE-chain} of a
$0$-$1$-filling is a sequence of $1$'s in the filling such that any $1$
in the sequence is above and to the right of the preceding $1$ in 
the sequence. Similarly, a {\it SE-chain} of a
$0$-$1$-filling is a set of $1$'s in the filling such that any $1$
in the sequence is below and to the right of the preceding $1$ in the sequence.

\proclaim{Theorem \TB} 
Given a growth diagram on a cell arrangement 
with empty partitions labelling all the
corners along the left boundary and the bottom boundary of the cell arrangement,
the partition $\la=(\la_1,\la_2,\dots,\la_\ell)$ labelling corner $c$ 
satisfies the following two properties:
\roster
\item"(G1)"For any $k$, the maximal cardinality of the union of $k$
NE-chains situated in the rectangular region to the left and
below of $c$ is equal to $\la_1+\la_2+\dots+\la_k$.
\item"(G2)"For any $k$, the maximal cardinality of the union of $k$
SE-chains situated in the rectangular region to the left and
below of $c$ is equal to $\la'_1+\la'_2+\dots+\la'_k$, where $\la'$
denotes the partition conjugate to $\la$.
\endroster
In particular, $\la_1$ is the length of the longest NE-chain
in the rectangular region to the left and below of $c$, and $\la'_1$
is the length of the longest SE-chain in the same
rectangular region.
\endproclaim

\subhead 3. The main theorems \endsubhead
Here, we state and prove our main results.
The theorem below proves and, at the same time, refines Burrill's
conjecture from the introduction. In the statement of the theorem and later,
the symbol $(1^m)$ stands for the partition 
$(1,1,\dots,1)$, with $m$ components
$1$, that is, the one-column shape of length~$m$.

\proclaim{Theorem \TC}
Let $n,m,k$ be given non-negative integers. 
The number of oscillating tableaux
of length $n$ with at most $k$ columns,
starting at $\emptyset$ and ending at the one-column shape $(1^m)$,
is equal to the number of standard Young tableaux of size~$n$
with $m$ columns of odd length, all columns of length at most $2k$.
\endproclaim

\demo{Proof}
We start with a standard Young tableau $T$ of size $n$ with at most $2k$ rows
and with $m$ columns of odd length. As a running example, we choose
$$
\matrix 1&3&4&8\\
2&6&7\\
5&10\\
9&12\\
11\endmatrix
$$
with $n=12$, $k=3$, and $m=2$. Indeed, this standard Young tableau has
$12$ entries, it has less than
or equal to $2k=6$ rows (namely~$5$) and $2$ columns of odd length.

\medskip
{\smc Step 1}. At the end of the odd-length columns, 
we put $I$, $II$, $III$, \dots,
from left to right. In our running example, we obtain
$$
\matrix 1&3&4&8\\
2&6&7&II\\
5&10\\
9&12\\
11\\
I\endmatrix
$$
We slide $I$, $II$, $III$, \dots, in this order, to the top-left of
the tableau, using (inverse) jeu de taquin (cf.\ \cite{\SagaAQ,
Sec.~3.7}).
That is, as long as above $I$
we find an entry belonging to $\{1,2,\dots,n\}$, we interchange $I$
with this entry; then, as long as to the left of or above $II$
we find entries belonging to $\{1,2,\dots,n\}$, we interchange $II$
with the larger of the two entries; then we do the same with $III$,
$IV$, etc. In the end, we obtain the
standard tableau $T'$ in the alphabet $I,II,III,\dots,1,2,\dots$.
In our running example, we get
$$
\matrix I&II&3&4\\
1&6&7&8\\
2&10\\
5&12\\
9\\
11\endmatrix
$$
It should be observed that, necessarily, $I$, $II$, $III$, \dots\
come to rest in the first row, in this order. The reason is
that successive jeu de taquin paths cannot cross each other, of which
one can easily convince oneself. More precisely, 
the jeu de taquin path of $I$ has to stay (weakly) to the left of the path
of~$II$, but strictly to the left along vertical pieces, the same
being true for the jeu de taquin paths of $II$ and $III$, etc.

\medskip
{\smc Step 2}. We interpret the tableau $T'$ as a sequence of
partitions 
$\emptyset=\la^0\subseteq\la^1\subseteq\dots\subseteq\la^{n+m}=\la$,
as described in Section~\SB. The partitions are placed along the
corners of the upper boundary of an $(n+m)\times(n+m)$ square
cell arrangement from right to left, 
and along the left boundary from bottom to top.
Then we apply the inverse growth diagram algorithm described in
Section~\SB, however, not in direction bottom/left but instead in direction
bottom/right. For the result in our running example see Figure~\FD.
\footnote{Alternatively, we could have applied
the inverse Robinson--Schensted algorithm to
the tableau pair $(T',T')$. This yields an involution
(since we started with a tableau pair consisting of two identical
tableaux) on $I,II,III,\dots,1,2,\dots$.
In our running example, we obtain the involution
$$
(I,5)(II,11)(1,9)(2,6)(3,7)(4,12)(8,10).
$$
We represent this involution as a 0-1-filling of a square,
where we label the rows $I,II,III,\dots,1,2,\dots$ from bottom to
top, and columns by the same labels from right to left
(skipping unconcerned labels). In our running example, this leads
to the $0$-$1$-filling in Figure~\FD\ (with $X$'s representing the $1$'s, while
empty cells represent the $0$'s). That this is indeed equivalent is
due to the fact (cf.\ \cite{\BrFoAA, pp.~95--98},
\cite{\StanBI, Theorem~7.13.5})
that the bijection between permutations and pairs of standard tableaux
defined by the growth diagram on the square coincides with
the Robinson--Schensted correspondence.}
In the figure, we have ``separated" the rows and columns corresponding
to the ``extra" letters $I$ and $II$ by thick lines.

\midinsert
$$
\Einheit.3cm
\Pfad(0,0),111111111111111111111111111111111111111111\endPfad
\Pfad(0,3),111111111111111111111111111111111111111111\endPfad
\Pfad(0,9),111111111111111111111111111111111111111111\endPfad
\Pfad(0,12),111111111111111111111111111111111111111111\endPfad
\Pfad(0,15),111111111111111111111111111111111111111111\endPfad
\Pfad(0,18),111111111111111111111111111111111111111111\endPfad
\Pfad(0,21),111111111111111111111111111111111111111111\endPfad
\Pfad(0,24),111111111111111111111111111111111111111111\endPfad
\Pfad(0,27),111111111111111111111111111111111111111111\endPfad
\Pfad(0,30),111111111111111111111111111111111111111111\endPfad
\Pfad(0,33),111111111111111111111111111111111111111111\endPfad
\Pfad(0,36),111111111111111111111111111111111111111111\endPfad
\Pfad(0,39),111111111111111111111111111111111111111111\endPfad
\Pfad(0,42),111111111111111111111111111111111111111111\endPfad
\Pfad(0,0),222222222222222222222222222222222222222222\endPfad
\Pfad(3,0),222222222222222222222222222222222222222222\endPfad
\Pfad(6,0),222222222222222222222222222222222222222222\endPfad
\Pfad(9,0),222222222222222222222222222222222222222222\endPfad
\Pfad(12,0),222222222222222222222222222222222222222222\endPfad
\Pfad(15,0),222222222222222222222222222222222222222222\endPfad
\Pfad(18,0),222222222222222222222222222222222222222222\endPfad
\Pfad(21,0),222222222222222222222222222222222222222222\endPfad
\Pfad(24,0),222222222222222222222222222222222222222222\endPfad
\Pfad(27,0),222222222222222222222222222222222222222222\endPfad
\Pfad(30,0),222222222222222222222222222222222222222222\endPfad
\Pfad(33,0),222222222222222222222222222222222222222222\endPfad
\Pfad(39,0),222222222222222222222222222222222222222222\endPfad
\Pfad(42,0),222222222222222222222222222222222222222222\endPfad
\PfadDicke{2pt}
\Pfad(0,6),111111111111111111111111111111111111111111\endPfad
\Pfad(36,0),222222222222222222222222222222222222222222\endPfad
\Label\ro{\text {\seventeenpoint X}}(22,1)
\Label\ro{\text {\seventeenpoint X}}(4,4)
\Label\ro{\text {\seventeenpoint X}}(10,7)
\Label\ro{\text {\seventeenpoint X}}(19,10)
\Label\ro{\text {\seventeenpoint X}}(16,13)
\Label\ro{\text {\seventeenpoint X}}(1,16)
\Label\ro{\text {\seventeenpoint X}}(40,19)
\Label\ro{\text {\seventeenpoint X}}(31,22)
\Label\ro{\text {\seventeenpoint X}}(28,25)
\Label\ro{\text {\seventeenpoint X}}(7,28)
\Label\ro{\text {\seventeenpoint X}}(34,31)
\Label\ro{\text {\seventeenpoint X}}(13,34)
\Label\ro{\text {\seventeenpoint X}}(37,37)
\Label\ro{\text {\seventeenpoint X}}(25,40)
\Label\u{\eightpoint\emptyset}(0,0)
\Label\u{\eightpoint\emptyset}(3,0)
\Label\u{\eightpoint\emptyset}(6,0)
\Label\u{\eightpoint\emptyset}(9,0)
\Label\u{\eightpoint\emptyset}(12,0)
\Label\u{\eightpoint\emptyset}(15,0)
\Label\u{\eightpoint\emptyset}(18,0)
\Label\u{\eightpoint\emptyset}(21,0)
\Label\u{\eightpoint\emptyset}(24,0)
\Label\u{\eightpoint\emptyset}(27,0)
\Label\u{\eightpoint\emptyset}(30,0)
\Label\u{\eightpoint\emptyset}(33,0)
\Label\u{\eightpoint\emptyset}(36,0)
\Label\u{\eightpoint\emptyset}(39,0)
\Label\ru{\eightpoint\emptyset}(42,0)
\Label\r{\eightpoint\emptyset}(42,3)
\Label\r{\eightpoint\emptyset}(42,6)
\Label\r{\eightpoint\emptyset}(42,9)
\Label\r{\eightpoint\emptyset}(42,12)
\Label\r{\eightpoint\emptyset}(42,15)
\Label\r{\eightpoint\emptyset}(42,18)
\Label\r{\eightpoint\emptyset}(42,21)
\Label\r{\eightpoint\emptyset}(42,24)
\Label\r{\eightpoint\emptyset}(42,27)
\Label\r{\eightpoint\emptyset}(42,30)
\Label\r{\eightpoint\emptyset}(42,33)
\Label\r{\eightpoint\emptyset}(42,36)
\Label\r{\eightpoint\emptyset}(42,39)
\Label\r{\eightpoint\emptyset}(42,42)
\Label\lo{\eightpoint\emptyset}(39,3)
\Label\lo{\eightpoint\emptyset}(36,3)
\Label\lo{\eightpoint\emptyset}(33,3)
\Label\lo{\eightpoint\emptyset}(30,3)
\Label\lo{\eightpoint\emptyset}(27,3)
\Label\lo{\eightpoint\emptyset}(24,3)
\Label\lo{\eightpoint1}(21,3)
\Label\lo{\eightpoint1}(18,3)
\Label\lo{\eightpoint1}(15,3)
\Label\lo{\eightpoint1}(12,3)
\Label\lo{\eightpoint1}(9,3)
\Label\lo{\eightpoint1}(6,3)
\Label\lo{\eightpoint1}(3,3)
\Label\lo{\eightpoint1}(0,3)
\Label\lo{\eightpoint\emptyset}(39,6)
\Label\lo{\eightpoint\emptyset}(36,6)
\Label\lo{\eightpoint\emptyset}(33,6)
\Label\lo{\eightpoint\emptyset}(30,6)
\Label\lo{\eightpoint\emptyset}(27,6)
\Label\lo{\eightpoint\emptyset}(24,6)
\Label\lo{\eightpoint1}(21,6)
\Label\lo{\eightpoint1}(18,6)
\Label\lo{\eightpoint1}(15,6)
\Label\lo{\eightpoint1}(12,6)
\Label\lo{\eightpoint1}(9,6)
\Label\lo{\eightpoint1}(6,6)
\Label\lo{\eightpoint2}(3,6)
\Label\lo{\eightpoint2}(0,6)
\Label\lo{\eightpoint\emptyset}(39,9)
\Label\lo{\eightpoint\emptyset}(36,9)
\Label\lo{\eightpoint\emptyset}(33,9)
\Label\lo{\eightpoint\emptyset}(30,9)
\Label\lo{\eightpoint\emptyset}(27,9)
\Label\lo{\eightpoint\emptyset}(24,9)
\Label\lo{\eightpoint1}(21,9)
\Label\lo{\eightpoint1}(18,9)
\Label\lo{\eightpoint1}(15,9)
\Label\lo{\eightpoint1}(12,9)
\Label\lo{\eightpoint2}(9,9)
\Label\lo{\eightpoint2}(6,9)
\Label\lo{\eightpoint21\ }(3,9)
\Label\lo{\eightpoint21\ }(0,9)
\Label\lo{\eightpoint\emptyset}(39,12)
\Label\lo{\eightpoint\emptyset}(36,12)
\Label\lo{\eightpoint\emptyset}(33,12)
\Label\lo{\eightpoint\emptyset}(30,12)
\Label\lo{\eightpoint\emptyset}(27,12)
\Label\lo{\eightpoint\emptyset}(24,12)
\Label\lo{\eightpoint1}(21,12)
\Label\lo{\eightpoint2}(18,12)
\Label\lo{\eightpoint2}(15,12)
\Label\lo{\eightpoint2}(12,12)
\Label\lo{\eightpoint21\ }(9,12)
\Label\lo{\eightpoint21\ }(6,12)
\Label\lo{\eightpoint211\ \ }(3,12)
\Label\lo{\eightpoint211\ \ }(0,12)
\Label\lo{\eightpoint\emptyset}(39,15)
\Label\lo{\eightpoint\emptyset}(36,15)
\Label\lo{\eightpoint\emptyset}(33,15)
\Label\lo{\eightpoint\emptyset}(30,15)
\Label\lo{\eightpoint\emptyset}(27,15)
\Label\lo{\eightpoint\emptyset}(24,15)
\Label\lo{\eightpoint1}(21,15)
\Label\lo{\eightpoint2}(18,15)
\Label\lo{\eightpoint3}(15,15)
\Label\lo{\eightpoint3}(12,15)
\Label\lo{\eightpoint31\ }(9,15)
\Label\lo{\eightpoint31\ }(6,15)
\Label\lo{\eightpoint311\ \ }(3,15)
\Label\lo{\eightpoint311\ \ }(0,15)
\Label\lo{\eightpoint\emptyset}(39,18)
\Label\lo{\eightpoint\emptyset}(36,18)
\Label\lo{\eightpoint\emptyset}(33,18)
\Label\lo{\eightpoint\emptyset}(30,18)
\Label\lo{\eightpoint\emptyset}(27,18)
\Label\lo{\eightpoint\emptyset}(24,18)
\Label\lo{\eightpoint1}(21,18)
\Label\lo{\eightpoint2}(18,18)
\Label\lo{\eightpoint3}(15,18)
\Label\lo{\eightpoint3}(12,18)
\Label\lo{\eightpoint31\ }(9,18)
\Label\lo{\eightpoint31\ }(6,18)
\Label\lo{\eightpoint311\ \ }(3,18)
\Label\lo{\eightpoint411\ \ }(0,18)
\Label\lo{\eightpoint1}(39,21)
\Label\lo{\eightpoint1}(36,21)
\Label\lo{\eightpoint1}(33,21)
\Label\lo{\eightpoint1}(30,21)
\Label\lo{\eightpoint1}(27,21)
\Label\lo{\eightpoint1}(24,21)
\Label\lo{\eightpoint11\ }(21,21)
\Label\lo{\eightpoint21\ }(18,21)
\Label\lo{\eightpoint31\ }(15,21)
\Label\lo{\eightpoint31\ }(12,21)
\Label\lo{\eightpoint311\ \ }(9,21)
\Label\lo{\eightpoint311\ \ }(6,21)
\Label\lo{\eightpoint3111\ \,\ \ }(3,21)
\Label\lo{\eightpoint4111\ \,\ \ }(0,21)
\Label\lo{\eightpoint1}(39,24)
\Label\lo{\eightpoint1}(36,24)
\Label\lo{\eightpoint1}(33,24)
\Label\lo{\eightpoint2}(30,24)
\Label\lo{\eightpoint2}(27,24)
\Label\lo{\eightpoint2}(24,24)
\Label\lo{\eightpoint21\ }(21,24)
\Label\lo{\eightpoint22\ }(18,24)
\Label\lo{\eightpoint32\ }(15,24)
\Label\lo{\eightpoint32\ }(12,24)
\Label\lo{\eightpoint321\ \ }(9,24)
\Label\lo{\eightpoint321\ \ }(6,24)
\Label\lo{\eightpoint3211\ \,\ \ }(3,24)
\Label\lo{\eightpoint4211\ \,\ \ }(0,24)
\Label\lo{\eightpoint1}(39,27)
\Label\lo{\eightpoint1}(36,27)
\Label\lo{\eightpoint1}(33,27)
\Label\lo{\eightpoint2}(30,27)
\Label\lo{\eightpoint3}(27,27)
\Label\lo{\eightpoint3}(24,27)
\Label\lo{\eightpoint31\ }(21,27)
\Label\lo{\eightpoint32\ }(18,27)
\Label\lo{\eightpoint33\ }(15,27)
\Label\lo{\eightpoint33\ }(12,27)
\Label\lo{\eightpoint331\ \ }(9,27)
\Label\lo{\eightpoint331\ \ }(6,27)
\Label\lo{\eightpoint3311\ \,\ \ }(3,27)
\Label\lo{\eightpoint4311\ \,\ \ }(0,27)
\Label\lo{\eightpoint1}(39,30)
\Label\lo{\eightpoint1}(36,30)
\Label\lo{\eightpoint1}(33,30)
\Label\lo{\eightpoint2}(30,30)
\Label\lo{\eightpoint3}(27,30)
\Label\lo{\eightpoint3}(24,30)
\Label\lo{\eightpoint31\ }(21,30)
\Label\lo{\eightpoint32\ }(18,30)
\Label\lo{\eightpoint33\ }(15,30)
\Label\lo{\eightpoint33\ }(12,30)
\Label\lo{\eightpoint331\ \ }(9,30)
\Label\lo{\eightpoint431\ \ }(6,30)
\Label\lo{\eightpoint4311\ \,\ \ }(3,30)
\Label\lo{\eightpoint4411\ \,\ \ }(0,30)
\Label\lo{\eightpoint1}(39,33)
\Label\lo{\eightpoint1}(36,33)
\Label\lo{\eightpoint2}(33,33)
\Label\lo{\eightpoint21\ }(30,33)
\Label\lo{\eightpoint31\ }(27,33)
\Label\lo{\eightpoint31\ }(24,33)
\Label\lo{\eightpoint311\ \ }(21,33)
\Label\lo{\eightpoint321\ \ }(18,33)
\Label\lo{\eightpoint331\ \ }(15,33)
\Label\lo{\eightpoint331\ \ }(12,33)
\Label\lo{\eightpoint3311\ \,\ \ }(9,33)
\Label\lo{\eightpoint4311\ \,\ \ }(6,33)
\Label\lo{\eightpoint43111\ \ \,\ \ }(3,33)
\Label\lo{\eightpoint44111\ \ \,\ \ }(0,33)
\Label\lo{\eightpoint1}(39,36)
\Label\lo{\eightpoint1}(36,36)
\Label\lo{\eightpoint2}(33,36)
\Label\lo{\eightpoint21\ }(30,36)
\Label\lo{\eightpoint31\ }(27,36)
\Label\lo{\eightpoint31\ }(24,36)
\Label\lo{\eightpoint311\ \ }(21,36)
\Label\lo{\eightpoint321\ \ }(18,36)
\Label\lo{\eightpoint331\ \ }(15,36)
\Label\lo{\eightpoint431\ \ }(12,36)
\Label\lo{\eightpoint4311\ \,\ \ }(9,36)
\Label\lo{\eightpoint4411\ \,\ \ }(6,36)
\Label\lo{\eightpoint44111\ \ \,\ \ }(3,36)
\Label\lo{\eightpoint44211\ \ \,\ \ }(0,36)
\Label\lo{\eightpoint1}(39,39)
\Label\lo{\eightpoint2}(36,39)
\Label\lo{\eightpoint21\ }(33,39)
\Label\lo{\eightpoint211\ \ }(30,39)
\Label\lo{\eightpoint311\ \ }(27,39)
\Label\lo{\eightpoint311\ \ }(24,39)
\Label\lo{\eightpoint3111\ \,\ \ }(21,39)
\Label\lo{\eightpoint3211\ \,\ \ }(18,39)
\Label\lo{\eightpoint3311\ \,\ \ }(15,39)
\Label\lo{\eightpoint4311\ \,\ \ }(12,39)
\Label\lo{\eightpoint43111\ \ \,\ \ }(9,39)
\Label\lo{\eightpoint44111\ \ \,\ \ }(6,39)
\Label\lo{\eightpoint441111\ \,\ \ \,\ \ }(3,39)
\Label\lo{\eightpoint442111\ \,\ \ \,\ \ }(0,39)
\Label\lo{\eightpoint1}(39,42)
\Label\lo{\eightpoint2}(36,42)
\Label\lo{\eightpoint21\ }(33,42)
\Label\lo{\eightpoint211\ \ }(30,42)
\Label\lo{\eightpoint311\ \ }(27,42)
\Label\lo{\eightpoint411\ \ }(24,42)
\Label\lo{\eightpoint4111\ \,\ \ }(21,42)
\Label\lo{\eightpoint4211\ \,\ \ }(18,42)
\Label\lo{\eightpoint4311\ \,\ \ }(15,42)
\Label\lo{\eightpoint4411\ \,\ \ }(12,42)
\Label\lo{\eightpoint44111\ \ \,\ \ }(9,42)
\Label\lo{\eightpoint44211\ \ \,\ \ }(6,42)
\Label\lo{\eightpoint442111\ \,\ \ \,\ \ }(3,42)
\Label\lo{\eightpoint442211\ \,\ \ \,\ \ }(0,42)
\hskip12.3cm
$$
\vskip10pt
\centerline{\smc Figure \FD}
\endinsert

Clearly, the growth diagram is symmetric with respect
to the main diagonal (i.e., the top/left--bottom/right diagonal). Furthermore,
by Greene's theorem, the lengths of north-east chains of $1$'s (i.e.,
of X's) are at most $2k$. 
Moreover, again by Greene's theorem (but now using item~(G2)),
since we started with a tableau with all columns of even length, 
there cannot be any crosses along the main diagonal.
(This is a fact observed earlier by Sch\"utzenberger
\cite{\SchuAB, p.~127} in a more general form.)

Finally, since the letters $I,II,III,\dots$ appeared in the
first row of $T'$ in that order, the $X$'s in the rows labelled by these
letters (the last rows; in the running example these are the last two rows)
form {\it one} SE-chain. (Also this follows from Greene's theorem.)
Together with the earlier observations that the diagram is symmetric
with respect to the main diagonal and that there are no crosses along
the main diagonal, it follows that the region
below the thick horizontal line and to the right of the thick vertical
horizontal line does not contain any crosses.

\medskip
{\smc Step 3}. Since the diagram is symmetric without crosses on the main
diagonal, we may forget about one half of the diagram, say the upper
half (including the main diagonal). 
In our running example, we arrive at the filling
of the staircase diagram displayed in
Figure~\FE.

\midinsert
$$
\Einheit.3cm
\Pfad(0,0),111111111111111111111111111111111111111111\endPfad
\Pfad(0,3),111111111111111111111111111111111111111\endPfad
\Pfad(0,9),111111111111111111111111111111111\endPfad
\Pfad(0,12),111111111111111111111111111111\endPfad
\Pfad(0,15),111111111111111111111111111\endPfad
\Pfad(0,18),111111111111111111111111\endPfad
\Pfad(0,21),111111111111111111111\endPfad
\Pfad(0,24),111111111111111111\endPfad
\Pfad(0,27),111111111111111\endPfad
\Pfad(0,30),111111111111\endPfad
\Pfad(0,33),111111111\endPfad
\Pfad(0,36),111111\endPfad
\Pfad(0,39),111\endPfad
\Pfad(0,0),222222222222222222222222222222222222222222\endPfad
\Pfad(3,0),222222222222222222222222222222222222222\endPfad
\Pfad(6,0),222222222222222222222222222222222222\endPfad
\Pfad(9,0),222222222222222222222222222222222\endPfad
\Pfad(12,0),222222222222222222222222222222\endPfad
\Pfad(15,0),222222222222222222222222222\endPfad
\Pfad(18,0),222222222222222222222222\endPfad
\Pfad(21,0),222222222222222222222\endPfad
\Pfad(24,0),222222222222222222\endPfad
\Pfad(27,0),222222222222222\endPfad
\Pfad(30,0),222222222222\endPfad
\Pfad(33,0),222222222\endPfad
\Pfad(39,0),222\endPfad
\PfadDicke{2pt}
\Pfad(0,6),111111111111111111111111111111111111\endPfad
\Pfad(36,0),222222\endPfad
\Label\ro{\text {\seventeenpoint X}}(22,1)
\Label\ro{\text {\seventeenpoint X}}(4,4)
\Label\ro{\text {\seventeenpoint X}}(10,7)
\Label\ro{\text {\seventeenpoint X}}(19,10)
\Label\ro{\text {\seventeenpoint X}}(16,13)
\Label\ro{\text {\seventeenpoint X}}(1,16)
\Label\ro{\text {\seventeenpoint X}}(7,28)
\hskip12.3cm
$$
\vskip10pt
\centerline{\smc Figure \FE}
\endinsert

Now we place empty partitions along the corners of the left and
the bottom boundary of the staircase diagram.
To the resulting diagram we apply the (forward) growth diagram construction
as described in Section~\SB, here in direction top/right.
Figure~\FF\ shows the result in our running example.

\vskip10pt
$$
\Einheit.3cm
\Pfad(0,0),111111111111111111111111111111111111111111\endPfad
\Pfad(0,3),111111111111111111111111111111111111111\endPfad
\Pfad(0,9),111111111111111111111111111111111\endPfad
\Pfad(0,12),111111111111111111111111111111\endPfad
\Pfad(0,15),111111111111111111111111111\endPfad
\Pfad(0,18),111111111111111111111111\endPfad
\Pfad(0,21),111111111111111111111\endPfad
\Pfad(0,24),111111111111111111\endPfad
\Pfad(0,27),111111111111111\endPfad
\Pfad(0,30),111111111111\endPfad
\Pfad(0,33),111111111\endPfad
\Pfad(0,36),111111\endPfad
\Pfad(0,39),111\endPfad
\Pfad(0,0),222222222222222222222222222222222222222222\endPfad
\Pfad(3,0),222222222222222222222222222222222222222\endPfad
\Pfad(6,0),222222222222222222222222222222222222\endPfad
\Pfad(9,0),222222222222222222222222222222222\endPfad
\Pfad(12,0),222222222222222222222222222222\endPfad
\Pfad(15,0),222222222222222222222222222\endPfad
\Pfad(18,0),222222222222222222222222\endPfad
\Pfad(21,0),222222222222222222222\endPfad
\Pfad(24,0),222222222222222222\endPfad
\Pfad(27,0),222222222222222\endPfad
\Pfad(30,0),222222222222\endPfad
\Pfad(33,0),222222222\endPfad
\Pfad(39,0),222\endPfad
\PfadDicke{2pt}
\Pfad(0,6),111111111111111111111111111111111111\endPfad
\Pfad(36,0),222222\endPfad
\Label\ro{\text {\seventeenpoint X}}(22,1)
\Label\ro{\text {\seventeenpoint X}}(4,4)
\Label\ro{\text {\seventeenpoint X}}(10,7)
\Label\ro{\text {\seventeenpoint X}}(19,10)
\Label\ro{\text {\seventeenpoint X}}(16,13)
\Label\ro{\text {\seventeenpoint X}}(1,16)
\Label\ro{\text {\seventeenpoint X}}(7,28)
\Label\u{\eightpoint\emptyset}(0,0)
\Label\u{\eightpoint\emptyset}(3,0)
\Label\u{\eightpoint\emptyset}(6,0)
\Label\u{\eightpoint\emptyset}(9,0)
\Label\u{\eightpoint\emptyset}(12,0)
\Label\u{\eightpoint\emptyset}(15,0)
\Label\u{\eightpoint\emptyset}(18,0)
\Label\u{\eightpoint\emptyset}(21,0)
\Label\u{\eightpoint\emptyset}(24,0)
\Label\u{\eightpoint\emptyset}(27,0)
\Label\u{\eightpoint\emptyset}(30,0)
\Label\u{\eightpoint\emptyset}(33,0)
\Label\u{\eightpoint\emptyset}(36,0)
\Label\u{\eightpoint\emptyset}(39,0)
\Label\u{\eightpoint\emptyset}(42,0)
\Label\l{\eightpoint\emptyset}(0,3)
\Label\l{\eightpoint\emptyset}(0,6)
\Label\l{\eightpoint\emptyset}(0,9)
\Label\l{\eightpoint\emptyset}(0,12)
\Label\l{\eightpoint\emptyset}(0,15)
\Label\l{\eightpoint\emptyset}(0,18)
\Label\l{\eightpoint\emptyset}(0,21)
\Label\l{\eightpoint\emptyset}(0,24)
\Label\l{\eightpoint\emptyset}(0,27)
\Label\l{\eightpoint\emptyset}(0,30)
\Label\l{\eightpoint\emptyset}(0,33)
\Label\l{\eightpoint\emptyset}(0,36)
\Label\l{\eightpoint\emptyset}(0,39)
\Label\l{\eightpoint\emptyset}(0,42)
\Label\ro{\eightpoint\emptyset}(3,3)
\Label\ro{\eightpoint\emptyset}(6,3)
\Label\ro{\eightpoint\emptyset}(9,3)
\Label\ro{\eightpoint\emptyset}(12,3)
\Label\ro{\eightpoint\emptyset}(15,3)
\Label\ro{\eightpoint\emptyset}(18,3)
\Label\ro{\eightpoint\emptyset}(21,3)
\Label\ro{\eightpoint1}(24,3)
\Label\ro{\eightpoint1}(27,3)
\Label\ro{\eightpoint1}(30,3)
\Label\ro{\eightpoint1}(33,3)
\Label\ro{\eightpoint1}(36,3)
\Label\ro{\eightpoint\emptyset}(3,6)
\Label\ro{\eightpoint1}(6,6)
\Label\ro{\eightpoint1}(9,6)
\Label\ro{\eightpoint1}(12,6)
\Label\ro{\eightpoint1}(15,6)
\Label\ro{\eightpoint1}(18,6)
\Label\ro{\eightpoint1}(21,6)
\Label\ro{\eightpoint11}(24,6)
\Label\ro{\eightpoint11}(27,6)
\Label\ro{\eightpoint11}(30,6)
\Label\ro{\eightpoint11}(33,6)
\Label\ro{\eightpoint\emptyset}(3,9)
\Label\ro{\eightpoint1}(6,9)
\Label\ro{\eightpoint1}(9,9)
\Label\ro{\eightpoint2}(12,9)
\Label\ro{\eightpoint2}(15,9)
\Label\ro{\eightpoint2}(18,9)
\Label\ro{\eightpoint2}(21,9)
\Label\ro{\eightpoint21}(24,9)
\Label\ro{\eightpoint21}(27,9)
\Label\ro{\eightpoint21}(30,9)
\Label\ro{\eightpoint\emptyset}(3,12)
\Label\ro{\eightpoint1}(6,12)
\Label\ro{\eightpoint1}(9,12)
\Label\ro{\eightpoint2}(12,12)
\Label\ro{\eightpoint2}(15,12)
\Label\ro{\eightpoint2}(18,12)
\Label\ro{\eightpoint3}(21,12)
\Label\ro{\eightpoint31}(24,12)
\Label\ro{\eightpoint31}(27,12)
\Label\ro{\eightpoint\emptyset}(3,15)
\Label\ro{\eightpoint1}(6,15)
\Label\ro{\eightpoint1}(9,15)
\Label\ro{\eightpoint2}(12,15)
\Label\ro{\eightpoint2}(15,15)
\Label\ro{\eightpoint3}(18,15)
\Label\ro{\eightpoint31}(21,15)
\Label\ro{\eightpoint\ \ 311}(24,15)
\Label\ro{\eightpoint1}(3,18)
\Label\ro{\eightpoint11}(6,18)
\Label\ro{\eightpoint11}(9,18)
\Label\ro{\eightpoint21}(12,18)
\Label\ro{\eightpoint21}(15,18)
\Label\ro{\eightpoint31}(18,18)
\Label\ro{\eightpoint\ \ 311}(21,18)
\Label\ro{\eightpoint1}(3,21)
\Label\ro{\eightpoint11}(6,21)
\Label\ro{\eightpoint11}(9,21)
\Label\ro{\eightpoint21}(12,21)
\Label\ro{\eightpoint21}(15,21)
\Label\ro{\eightpoint31}(18,21)
\Label\ro{\eightpoint1}(3,24)
\Label\ro{\eightpoint11}(6,24)
\Label\ro{\eightpoint11}(9,24)
\Label\ro{\eightpoint21}(12,24)
\Label\ro{\eightpoint21}(15,24)
\Label\ro{\eightpoint1}(3,27)
\Label\ro{\eightpoint11}(6,27)
\Label\ro{\eightpoint11}(9,27)
\Label\ro{\eightpoint21}(12,27)
\Label\ro{\eightpoint1}(3,30)
\Label\ro{\eightpoint11}(6,30)
\Label\ro{\eightpoint21}(9,30)
\Label\ro{\eightpoint1}(3,33)
\Label\ro{\eightpoint11}(6,33)
\Label\ro{\eightpoint1}(3,36)
\Label\ro{1}(3,39)
\Label\ro{11}(6,36)
\Label\ro{21}(9,33)
\Label\ro{22}(12,30)
\Label\ro{21}(15,27)
\Label\ro{31}(18,24)
\Label\ro{311}(21,21)
\Label\ro{3111}(24,18)
\Label\ro{311}(27,15)
\Label\ro{31}(30,12)
\Label\ro{21}(33,9)
\Label\ro{11}(36,6)
\Label\ro{1}(39,3)
\hskip13.1cm
$$
\vskip10pt
\centerline{\smc Figure \FF}
\vskip10pt

Along the corners on the main diagonal one reads an
oscillating tableau of length $n+m$. (In the running example in
Figure~\FF, this is the sequence of larger printed partitions.) 
However, since by one of the
previous observations (plus Greene's theorem) we know that the
last $m+1$ partitions (shapes) in the oscillating tableau will be
$(1^m),(1^{m-1}),(1^{m-2}),\dots,(1,1),(1),\emptyset$, we may discard all of
them except $(1^m)$, and in this way
obtain an oscillating tableau of length $n$, starting at $\emptyset$
and ending at $(1^m)$. Using Greene's theorem once more, we also see
that no shape along the main diagonal can have more than $k$ columns.

\medskip
Since every step in this construction can be reversed
in straightforward fashion, this yields the desired bijection.
\quad \quad \qed 
\enddemo

The announced ``Knuth-type" extension of Theorem~\TC\ is the following.

\proclaim{Theorem \TD}
Let $n,m,k$ and $j_1,j_2,\dots,j_n$ be non-negative integers.
The number of sequences of partitions
$\emptyset=\la^0,\la^1,\dots,\la^{2n}=(1^m)$
of length $2n$, with $\la^{2i-2}\supseteq\la^{2i-1}$ and
$\la^{2i-1}\subseteq\la^{2i}$ for $i=1,2,\dots,n$,
where each pair $(\la^{i-1},\la^i)$ differs by a vertical strip
{\rm(}that is, by a collection of cells which contains at most one  
cell in each row{\rm)}, $i=1,2,\dots,n$, 
where each partition $\la^i$ has at most $k$ columns,
and where 
$\vert\la^{2i-2}\vert-2\vert\la^{2i-1}\vert+\vert\la^{2i}\vert=j_{n-i+1}$
{\rm(}the left-hand side is the sum of the differences in sizes of
$(\la^{2i-2},\la^{2i-1})$ and of $(\la^{2i-1},\la^{2i})${\rm)},
$i=1,2,\dots,n$,
is equal to the number of semistandard tableaux with $j_i$ entries~$i$,
$i=1,2,\dots,n$,
with $m$ columns of odd length, all columns of length at most $2k$.
\endproclaim

\remark{Remark}
Theorem~\TC\ is the special case of Theorem~\TD\ where
$j_1=j_2=\dots=j_n=1$.
\endremark

\demo{Sketch of proof}
One proceeds in analogy with the proof of Theorem~\TC.
As a running example for illustration, we choose $n=4$, $m=2$, $k=2$, 
$j_1=5$, 
$j_2=2$, 
$j_3=6$, 
$j_4=3$, 
and the semistandard tableau
$$
\matrix 
1&1&1&1&1&3&3\\
2&2&3&3&4&4\\
3&3\\
4
\endmatrix
$$
Indeed, this semistandard tableau has $m=2$ columns of odd length,
all columns of length at most $2k=4$, 
$j_1=5$ entries~$1$, 
$j_2=2$ entries~$2$, 
$j_3=6$ entries~$3$, and
$j_4=3$ entries~$4$. 

\medskip
{\smc Step 1}. 
We place $I,II,\dots$ at the end of the columns of odd length,
from left to right. In our running example, we obtain
$$
\matrix 
1&1&1&1&1&3&3\\
2&2&3&3&4&4&II\\
3&3\\
4&I
\endmatrix
$$
Then we slide $I,II,\dots$ up to the first row. The result in our
example is
$$
\matrix 
I&II&1&1&1&3&3\\
1&1&2&3&3&4&4\\
2&3\\
3&4
\endmatrix
$$

\midinsert
$$
\Einheit.3cm
\Pfad(0,0),111111111111111111\endPfad
\Pfad(0,3),111111111111111111\endPfad
\Pfad(0,9),111111111111111111\endPfad
\Pfad(0,12),111111111111111111\endPfad
\Pfad(0,15),111111111111111111\endPfad
\Pfad(0,18),111111111111111111\endPfad
\Pfad(0,0),222222222222222222\endPfad
\Pfad(3,0),222222222222222222\endPfad
\Pfad(6,0),222222222222222222\endPfad
\Pfad(9,0),222222222222222222\endPfad
\Pfad(15,0),222222222222222222\endPfad
\Pfad(18,0),222222222222222222\endPfad
{\PfadDicke{2pt}
\Pfad(0,6),111111111111111111\endPfad
\Pfad(12,0),222222222222222222\endPfad}%
\Label\ro{\text {\seventeenpoint $0$}}(1,1)
\Label\ro{\text {\seventeenpoint $0$}}(4,1)
\Label\ro{\text {\seventeenpoint $0$}}(7,1)
\Label\ro{\text {\seventeenpoint $1$}}(10,1)
\Label\ro{\text {\seventeenpoint $0$}}(13,1)
\Label\ro{\text {\seventeenpoint $0$}}(16,1)
\Label\ro{\text {\seventeenpoint $0$}}(1,4)
\Label\ro{\text {\seventeenpoint $1$}}(4,4)
\Label\ro{\text {\seventeenpoint $0$}}(7,4)
\Label\ro{\text {\seventeenpoint $0$}}(10,4)
\Label\ro{\text {\seventeenpoint $0$}}(13,4)
\Label\ro{\text {\seventeenpoint $0$}}(16,4)
\Label\ro{\text {\seventeenpoint $1$}}(1,7)
\Label\ro{\text {\seventeenpoint $2$}}(4,7)
\Label\ro{\text {\seventeenpoint $1$}}(7,7)
\Label\ro{\text {\seventeenpoint $0$}}(10,7)
\Label\ro{\text {\seventeenpoint $0$}}(13,7)
\Label\ro{\text {\seventeenpoint $1$}}(16,7)
\Label\ro{\text {\seventeenpoint $0$}}(1,10)
\Label\ro{\text {\seventeenpoint $1$}}(4,10)
\Label\ro{\text {\seventeenpoint $0$}}(7,10)
\Label\ro{\text {\seventeenpoint $1$}}(10,10)
\Label\ro{\text {\seventeenpoint $0$}}(13,10)
\Label\ro{\text {\seventeenpoint $0$}}(16,10)
\Label\ro{\text {\seventeenpoint $2$}}(1,13)
\Label\ro{\text {\seventeenpoint $0$}}(4,13)
\Label\ro{\text {\seventeenpoint $1$}}(7,13)
\Label\ro{\text {\seventeenpoint $2$}}(10,13)
\Label\ro{\text {\seventeenpoint $1$}}(13,13)
\Label\ro{\text {\seventeenpoint $0$}}(16,13)
\Label\ro{\text {\seventeenpoint $0$}}(1,16)
\Label\ro{\text {\seventeenpoint $2$}}(4,16)
\Label\ro{\text {\seventeenpoint $0$}}(7,16)
\Label\ro{\text {\seventeenpoint $1$}}(10,16)
\Label\ro{\text {\seventeenpoint $0$}}(13,16)
\Label\ro{\text {\seventeenpoint $0$}}(16,16)
\Label\u{\eightpoint\emptyset}(0,0)
\Label\u{\eightpoint\emptyset}(3,0)
\Label\u{\eightpoint\emptyset}(6,0)
\Label\u{\eightpoint\emptyset}(9,0)
\Label\u{\eightpoint\emptyset}(12,0)
\Label\u{\eightpoint\emptyset}(15,0)
\Label\ru{\eightpoint\emptyset}(18,0)
\Label\r{\eightpoint\emptyset}(18,3)
\Label\r{\eightpoint\emptyset}(18,6)
\Label\r{\eightpoint\emptyset}(18,9)
\Label\r{\eightpoint\emptyset}(18,12)
\Label\r{\eightpoint\emptyset}(18,15)
\Label\r{\eightpoint\emptyset}(18,18)
\Label\lo{\eightpoint\emptyset}(15,3)
\Label\lo{\eightpoint\emptyset}(12,3)
\Label\lo{\eightpoint1}(9,3)
\Label\lo{\eightpoint1}(6,3)
\Label\lo{\eightpoint1}(3,3)
\Label\lo{\eightpoint1}(0,3)
\Label\lo{\eightpoint\emptyset}(15,6)
\Label\lo{\eightpoint\emptyset}(12,6)
\Label\lo{\eightpoint1}(9,6)
\Label\lo{\eightpoint1}(6,6)
\Label\lo{\eightpoint2}(3,6)
\Label\lo{\eightpoint2}(0,6)
\Label\lo{\eightpoint1}(15,9)
\Label\lo{\eightpoint1}(12,9)
\Label\lo{\eightpoint11\ }(9,9)
\Label\lo{\eightpoint21\ }(6,9)
\Label\lo{\eightpoint42\ }(3,9)
\Label\lo{\eightpoint52\ }(0,9)
\Label\lo{\eightpoint1}(15,12)
\Label\lo{\eightpoint1}(12,12)
\Label\lo{\eightpoint21\ }(9,12)
\Label\lo{\eightpoint22\ }(6,12)
\Label\lo{\eightpoint521\ \ }(3,12)
\Label\lo{\eightpoint531\ \ }(0,12)
\Label\lo{\eightpoint1}(15,15)
\Label\lo{\eightpoint2}(12,15)
\Label\lo{\eightpoint42\ }(9,15)
\Label\lo{\eightpoint521\ \ }(6,15)
\Label\lo{\eightpoint5511\ \,\ \ }(3,15)
\Label\lo{\eightpoint7521\ \,\ \ }(0,15)
\Label\lo{\eightpoint1}(15,18)
\Label\lo{\eightpoint2}(12,18)
\Label\lo{\eightpoint52\ }(9,18)
\Label\lo{\eightpoint531\ \ }(6,18)
\Label\lo{\eightpoint7521\ \,\ \ }(3,18)
\Label\lo{\eightpoint7722\ \,\ \ }(0,18)
\hbox{\hskip7.0cm}
\Pfad(0,0),111111111111111111\endPfad
\Pfad(0,3),111111111111111\endPfad
\Pfad(0,9),111111111\endPfad
\Pfad(0,12),111111\endPfad
\Pfad(0,15),111\endPfad
\Pfad(0,0),222222222222222222\endPfad
\Pfad(3,0),222222222222222\endPfad
\Pfad(6,0),222222222222\endPfad
\Pfad(9,0),222222222\endPfad
\Pfad(15,0),222\endPfad
{\PfadDicke{2pt}
\Pfad(0,6),111111111111\endPfad
\Pfad(12,0),222222\endPfad}%
\Label\ro{\text {\seventeenpoint $0$}}(1,1)
\Label\ro{\text {\seventeenpoint $0$}}(4,1)
\Label\ro{\text {\seventeenpoint $0$}}(7,1)
\Label\ro{\text {\seventeenpoint $1$}}(10,1)
\Label\ro{\text {\seventeenpoint $0$}}(13,1)
\Label\ro{\text {\seventeenpoint $0$}}(1,4)
\Label\ro{\text {\seventeenpoint $1$}}(4,4)
\Label\ro{\text {\seventeenpoint $0$}}(7,4)
\Label\ro{\text {\seventeenpoint $0$}}(10,4)
\Label\ro{\text {\seventeenpoint $1$}}(1,7)
\Label\ro{\text {\seventeenpoint $2$}}(4,7)
\Label\ro{\text {\seventeenpoint $1$}}(7,7)
\Label\ro{\text {\seventeenpoint $0$}}(1,10)
\Label\ro{\text {\seventeenpoint $1$}}(4,10)
\Label\ro{\text {\seventeenpoint $2$}}(1,13)
\Label\u{\eightpoint\emptyset}(0,0)
\Label\u{\eightpoint\emptyset}(3,0)
\Label\u{\eightpoint\emptyset}(6,0)
\Label\u{\eightpoint\emptyset}(9,0)
\Label\u{\eightpoint\emptyset}(12,0)
\Label\u{\eightpoint\emptyset}(15,0)
\Label\u{\eightpoint\emptyset}(18,0)
\Label\l{\eightpoint\emptyset}(0,3)
\Label\ro{\eightpoint\emptyset}(3,3)
\Label\ro{\eightpoint\emptyset}(6,3)
\Label\ro{\eightpoint\emptyset}(9,3)
\Label\ro{\eightpoint1}(12,3)
\Label\ro{\eightpoint1}(15,3)
\Label\l{\eightpoint\emptyset}(0,6)
\Label\ro{\eightpoint\emptyset}(3,6)
\Label\ro{\eightpoint1}(6,6)
\Label\ro{1}(9,6)
\Label\ro{\ 11}(12,6)
\Label\l{\eightpoint\emptyset}(0,9)
\Label\ro{\eightpoint1}(3,9)
\Label\ro{\ \ \ \ 1111}(6,9)
\Label\ro{\ \ \ \ 2111}(9,9)
\Label\l{\eightpoint\emptyset}(0,12)
\Label\ro{1}(3,12)
\Label\ro{\ \ \ \ 2111}(6,12)
\Label\l{\emptyset}(0,15)
\Label\ro{\ \ 111}(3,15)
\Label\l{\emptyset}(0,18)
\hskip5.2cm
$$
\vskip10pt
\centerline{\eightpoint a. Second step
\hskip5cm
b. Third step}
\vskip6pt
\centerline{\smc Figure \FG}
\endinsert

\medskip
{\smc Step 2}. 
Instead of applying the (ordinary) inverse growth diagram algorithm,
we apply its Knuth-type extension as described
in \cite{\KratCE, Sec.~4.1} or \cite{\RobyAA, Sec.~4.1}. 
Instead of a symmetric $0$-$1$-filling,
here we obtain a symmetric matrix with non-negative integer entries.
Again, along the diagonal there will be only $0$'s.
Figure~\FG.a shows what we
obtain in our running example.

\medskip
{\smc Step 3}. For the final step, in place of the
(ordinary) growth diagram algorithm, we apply its
Knuth-type extension described in \cite{\KratCE, Sec.~4.4} or
\cite{\LeeuAH, Sec.~3.2}.
The result for our running example is shown in
Figure~\FG.b. To obtain the (generalised) oscillating tableau, we
must read {\it all\/} partitions along the top-right boundary of
the cell arrangement (that is, including those which label ``inner"
corners), again skipping the ones in the region 
to the right of the thick vertical line. 
In Figure~\FG, these are the large printed
partitions, namely
$$
\emptyset,\ 
\emptyset,\
111,\
1,\
2111,\
1111,\
2111,\
1,\
11.
$$
It is not difficult to see that this map has all
the desired properties.
\quad \quad \qed
\enddemo

\subhead 4. Concluding remarks\endsubhead 

\medskip
(1) For both numbers in Theorem~\TC, there are explicit formulae
available. As pointed out in \cite{\BuCFAA, paragraphs around Theorem~19}, 
Gessel and Zeilberger's general result \cite{\GeZeAA} 
on enumeration of lattice paths by means of the reflection principle
yields (see Eq.~(38) in \cite{\GrMrAA}
with $n$ replaced by $k$, $\lambda=(k,k-1,\dots,1)$, and
$\eta=(m+k,k-1,k-2,\dots,1)$) that the number of oscillating
tableaux in Theorem~\TC\ is given by the coefficient of $t^n/n!$ in
$$
\det\left(
I_{i-j+m\cdot\chi(i=k)}(2t)
-
I_{i+j+m\cdot\chi(i=k)}(2t)
\right)_{1\le i,j\le k},
\tag\AA
$$
where $I_\alpha(x)$ is the modified Bessel function of the first kind
$$I_\al(x)=\sum _{\ell=0} ^{\infty}\frac {(x/2)^{2\ell+\al}}
{\ell!\,(\ell+\al)!},
$$
and $\chi(\Cal A)=1$ if $\Cal A$ is true and $\chi(\Cal A)=0$ otherwise.
On the other hand, one obtains the same formula for the number of
standard Young tableaux in Theorem~\TC\ from a result of
Goulden \cite{\GoulAD} (see \cite{\KratBC, Eq.~(3.6)} for a different
proof). Namely, if one extracts the coefficient of 
$x_1x_2\cdots x_n$ in Theorem~2.6 of 
\cite{\GoulAD} with $m$ and $k$ interchanged, and after having applied
the involution $\omega$ on symmetric functions which maps the
complete homogeneous symmetric functions to the elementary symmetric
functions (cf\. \cite{\MacdAC}), 
then the conclusion is that the number of standard Young
tableaux in Theorem~\TC\ is given by the coefficient of $t^n/n!$ in (\AA).
This amounts to an alternative --- albeit very roundabout and involved,
non-illuminating --- computational proof of Theorem~\TC. 
A proof of Theorem~\TD\ along the same lines is also possible.

\medskip
(2) What happens if we want to count standard Young tableaux
of size~$n$ with $m$ columns of odd length, all columns of length at
most $2k+1$ (instead of at most~$2k$)? As it turns out, the
corresponding number equals $\binom nm$ times 
the number of all standard Young tableaux
of size~$n-m$ with {\it all columns of even length not
exceeding~$2k$}. (In other words: the above problem can be reduced
to the $m=0$ case in Theorem~\TC.)
This is seen as follows: 
let $T$ be a standard Young tableau of size~$n$ with $m$ columns of
odd length, all columns of length at most $2k+1$.
To the last entries in the odd columns of $T$ one applies
the inverse mapping of Robinson--Schensted insertion
(cf\. \cite{\SagaAQ, Proof of Theorem~3.1.1}),
starting with the last entry in the right-most odd column, continuing
with the last entry in the (then) right-most odd column, until all
odd columns have disappeared. This produces a standard Young tableau
of size $n-m$ with only even columns, all of which have length at
most~$2k$, and a subset of $\{1,2,\dots,n\}$ of cardinality~$m$.
It is easy to see that all the steps in this mapping can be reversed
so that this describes a bijection. 

\medskip
(3) As indicated in the introduction, the difference between the
bijection proving Theorem~\TC\ presented here and that
of \cite{\BuCFAA} lies in the way one keeps track of 
the parameter~$m$ in Theorem~\TC.
Namely, in Step~1 of our proof of Theorem~\TC, we introduce the auxiliary
letters $I,II,III,\dots$ in order to ``make the $m$ odd columns even,"
and move the letters inside the tableau by jeu de taquin.
Then, in Step~2, we apply the (inverse) growth diagram construction to the
complete square, and finally, in Step~3, we apply the (forward) growth
diagram construction to half of the square to obtain the corresponding
oscillating tableau. 

On the other hand, if one realises the construction in \cite{\BuCFAA}
by means of growth diagrams, then
Burrill, Courtiel, Fusy, Melczer and Mishna apply the (inverse) growth
diagram construction directly, without any ``preprocessing."
The ``price to pay" is that they do obtain X's on the main diagonal.
These $m$ X's must be somehow moved ``into the half-square," and in order to
be able to do this without creating any unwanted chains, the growth
diagram construction has to be first played forth and back on the
half-square. Only then, the X's on the diagonal can be ``moved
inside," and a final application of the (forward) growth diagram
construction on the half-square completes the bijection. As far as I
can see, other than that, the two constructions do not seem to be more
deeply related.

\enddocument